%% file: main_arxiv.tex
\documentclass[a4paper,
10pt,
]{article}

\usepackage[english]{babel}

\usepackage[a4paper, left=2.5cm,right=2.5cm]{geometry}

\usepackage{hyperref}
\usepackage{bm}
\usepackage{authblk}
\usepackage{siunitx}
\usepackage{graphicx}
\usepackage[ruled]{algorithm2e}

\usepackage{caption}
\usepackage{subcaption}
\usepackage{amssymb} 
\usepackage[dvipsnames]{xcolor}
\usepackage{booktabs}
\usepackage{rotating}
\usepackage{tikz}
\usepackage[square,numbers]{natbib}
\usepackage{lineno}

\title{ \large \textbf{Non-intrusive model reduction of advection-dominated hyperbolic problems 
using neural network shift augmented manifold transformation}}

\author[1]{Harshith~Gowrachari\footnote{hgowrach@sissa.it}}
\author[1]{Nicola~Demo\footnote{ndemo@sissa.it}}
\author[2]{Giovanni~Stabile\footnote{giovanni.stabile@santannapisa.it}}
\author[1]{Gianluigi~Rozza\footnote{grozza@sissa.it}}

\affil[1]{\small Mathematics Area, mathLab, International School for Advanced Studies,
via Bonomea 265, 34136 Trieste, Italy}
\affil[2]{\small Biorobotics Institute, Sant’Anna School of Advanced Studies, V.le R. Piaggio 34, 56025, Pontedera, Pisa, Italy.}

\date{} 

\begin{document}
\maketitle
\input{sections/abstract}

\tableofcontents
\input{sections/intro}
\input{sections/methodology}

\input{sections/num_exp}
\input{sections/conclusion}
\input{sections/ackno}

\input{sections/appendix}

\bibliographystyle{abbrvnat}
\bibliography{bib/biblio}


\end{document}

%% file: sections/abstract.tex
\begin{abstract}

Advection-dominated problems are predominantly noticed in nature, engineering systems, and various industrial processes. Traditional linear compression methods, such as proper orthogonal decomposition (POD) and reduced basis (RB) methods are ill-suited for these problems, due to slow Kolmogorov \textit{n}-width decay. This results in inefficient and inaccurate reduced order models (ROMs). There are few non-linear approaches to accelerate the Kolmogorov \textit{n}-width decay. In this work, we use a neural network shift augmented transformation technique that employs automatic shift detection. \textcolor{Black}{This approach leverages a deep-learning framework to derive a parameter-dependent mapping between the original manifold $\mathcal{M}$ and the transformed manifold $\tilde\mathcal{M}$. We apply a linear compression method to obtain a low-dimensional linear approximation subspace of the transformed manifold $\tilde\mathcal{M}$. Furthermore, we construct non-intrusive reduced order models on the resulting transformed linear approximation subspace and employ automatic shift detection for predictions in the online stage. We propose a complete framework, the neural network shift-augmented proper orthogonal decomposition-based reduced order model (NNsPOD-ROM) algorithm, comprising both offline and online stages for model reduction of advection-dominated problems.} We test our proposed methodology on numerous experiments to evaluate its performance on the 1D linear advection equation, a higher order method benchmark case - the 2D isentropic convective vortex, and 2D two-phase flow.

\end{abstract}


%% file: sections/intro.tex
\section{Introduction}
Industries have a keen interest in obtaining reduced order models (ROMs) for engineering systems, particularly in control, optimization, and uncertainty quantification applications. ROMs provide computational efficiency in many query scenarios and are well-suited for real-time computations. \textcolor{Black}{These models create a low-dimensional representation of solutions to full-order models (FOMs), i.e., high-dimensional discretized partial differential equations (PDEs),  thereby significantly reducing the computational cost of numerical simulations compared to FOMs.}\\

To solve parametric PDEs that describe real-world scenarios, various discretization methods are used, such as the finite difference method (FDM), the finite element method (FEM), the finite volume method (FVM), the discontinuous Galerkin method (DGM), and the spectral element method (SEM). Although these methods are accurate, they are also computationally intensive. The ROM framework addresses this by dividing the process into two stages:

\begin{itemize}
    \item Offline stage (training): This expensive phase involves computing a set of FOM solutions (snapshots) and constructing the ROM based on a low-dimensional approximation subspace of the FOM solution manifold.
    \item Online stage (testing and predictive): This cost-effective phase uses the compressed information from the offline stage to predict reduced solutions for new unseen parameters, significantly reducing the computation costs.
\end{itemize}

There are various model reduction methods to obtain a low-dimensional approximation subspace \cite{quarteroni2015reduced, hesthaven2016certified, Chinesta2017, rozza2022advanced}, linear compression methods such as proper orthogonal decomposition (POD) and reduced basis (RB) methods are the most popular and promising techniques for various fluid dynamics and solid mechanics applications. Although linear compression methods are broadly effective across numerous applications, many problems of interest (e.g., advection-dominated problems) exhibit slow Kolmogorov \textit{n}-width decay, which limits the reduction attainable via linear approximation subspace \cite{Pinkus1985}.  \\

\textcolor{Black}{Advection-dominated problems exhibit slow Kolmogorov \textit{n}-width (KnW) decay. As a consequence, it is impossible to find a low-dimensional linear approximation subspace that yields accurate and efficient ROMs for these problems.} The KnW \cite{Pinkus1985} provides a rigorous measure of the reducibility of the solution manifold $\mathcal{M}$ by linear approximation subspaces. For the solution manifold $\mathcal{M}$ consisting of \textit{f} elements (FOM solutions) embedded in some normed linear space \( ({X}_{\mathcal{N}}, \| \cdot \|_{X_{\mathcal{N}}}) \), the  definition of Kolmogorov \textit{n}-width $d_n\left(\mathcal{M}\right)$ is as follows: 

\begin{equation}
    d_n\left(\mathcal{M}\right) := \inf _{E_n \subset X_N} \sup _{f \in \mathcal{M}_N} \inf _{g \in E_n}\|f-g\|_{X_{\mathcal{N}}},
\end{equation}

\noindent where ${E}_{n}$ is the linear subspace of $X_{\mathcal{N}}$ of dimension \textit{n}. The Kolmogorov \textit{n}-width $d_n(\mathcal{M})$ represents the worst error we commit by approximating any element \textit{f} of  $\mathcal{M}$ with approximated element \textit{g} of a linear subspace $E_n$ of dimension $n$. Moreover, the faster the KnW $d_n(\mathcal{M})$ decreases with increasing $n$, the low-dimensional linear approximation subspace of dimension $n$ becomes more effective at approximating the solution manifold $\mathcal{M}$. \textcolor{Black}{For the elliptic and parabolic problems, the KnW decay is exponential, i.e., $d_n(\mathcal{M}) \leq Ce^{-\beta n}$, where constants $C < \infty $ and $\beta > 0$ \cite{Buffa2012, ohlberger2016reduced}, which explains the effectiveness of model reduction methods for this class of problems. For hyperbolic problems, the n-width decay is at most $d_n(\mathcal{M}) > n ^{-1/2}$ \cite{Greif2019, ohlberger2016reduced}, which indicates that advection-dominated problems cannot be effectively approximated using low-dimensional linear subspaces.} Even the solution manifold of a simple linear advection problem exhibits slow KnW decay \cite{peherstorfer2022breaking}, meaning there is no low-dimensional linear approximation subspace that can effectively approximate the FOM solutions, and even if we consider the best-possible subspace, a large space should be used to obtain a small error. This poses a considerable challenge for effective model reduction, as computational efficiency depends on the reduced linear subspace. Consequently, linear compression methods such as POD and RB are ill-suited for this class of problems, resulting in inefficient and inaccurate linear projection-based ROMs. \\

In recent years there has been significant development in model-reduction methods to overcome the barrier of the slow KnW decay, which are leading towards non-linear approaches \cite{peherstorfer2022breaking}, particularly by performing model-reduction on non-linear manifolds, by obtaining non-linear trial subspaces using convolutional auto-encoders (CAEs) \cite{ kashima2016nonlinear, hartman2017deep, crisovan2019model, hoang2022projection, lee2020model, fresca2021comprehensive, kim2022fast, romor2023non} and by employing quadratic manifold \cite{Jain2017, Rutzmoser2017,  barnett2022quadratic, Geelen2023}; by incorporating adaptive local-in-parameter and/or local-in-space enrichment strategies \cite{haasdonk2008adaptive, washabaugh2012nonlinear, Carlberg2014, peherstorfer2015online, peherstorfer2020model, bruna2024neural} to derive efficient ROMs; by relying on multiple-linear subspaces (known as the dictionary approach) for constructing local approximation subspace, instead of single global approximation subspace \cite{Amsallem2012, Peherstorfer2014, Daniel2020, Geelen2022}; and by transforming the linear subspace to improve the its approximation properties for advection-dominated problems, such methods incorporate freezing \cite{ohlberger2013nonlinear}, shifting the POD basis \cite{reiss2018shifted} and manifold calibration/transformation \cite{mojgani2017arbitrary, iollo2014advection, cagniart2017model, rim2018transport, taddei2020registration, torlo2020model, nonino2023reduced, Nonino2024}. These approaches substantially rely on the prior knowledge of the underlying problem and physical properties such as advection phenomena that guide the basis shifting. \newline

\textcolor{Black}{In our work, we utilize the neural network shift augmented transformation technique employed within the neural network shifted-proper orthogonal decomposition (NNsPOD) framework \cite{papapicco2022neural}. This framework incorporates automatic shift detection to obtain a transformed linear subspace and does not require substantial prior knowledge of the underlying problem, such as governing equations or advection phenomena. This feature overcomes the limitation of the shifted-POD (sPOD) approach \cite{papapicco2022neural, reiss2018shifted}, which relies on prior knowledge of the problem like several other transformation approaches mentioned earlier. For a wide range of problems, knowing the shifts a priori is either challenging or impossible, where automatic shift detection plays a crucial role. The neural network shift augmented transformation technique seeks optimal mapping of numerous snapshots to reference frame by leveraging neural networks for automatic-shift detection, aiming to obtain low-dimensional transformed linear approximation subspace. Similarly, this type of transformation is achieved by solving an optimal transport problem as outlined in \cite{iollo2014advection, Iollo2022}, using a registration method \cite{taddei2020registration, Kleikamp2022} and through an implicit tracking method \cite{AlirezaMirhoseini2023}. In all these approaches, the mapping is built via an optimization procedure that aligns numerous snapshots with a reference state. Additionally, in \cite{Issan2023}, the method of characteristics and cross-correlation are utilised to circumvent the slow decay of KnW.} \newline

In reduced order modelling, once the low-dimensional approximation subspace of the FOM solution manifold is obtained, to compute the evolution of the dynamics in the low-dimensional subspace, modal coefficients or latent coordinates need to be computed. For this, there are intrusive and non-intrusive approaches. In the intrusive approach, the discretized PDEs (governing equations) are projected onto a low-dimensional subspace to obtain a system of low-dimensional ODEs, this procedure of ROM construction is known as Galerkin or Petrov-Galerkin projection-based ROM \cite{Benner2015, stabile2018finite, stabile2019reduced, lee2020model, romor2023non}. As this approach leverages equations, it is inherently physics-based. However, this method requires knowledge of the underlying equations and access to the discretized differential operators constructed by the full-order model, meaning complete access to the internal data structure of the solver, which makes it quite prohibitive for the application in industrial settings. For the case of non-linear/non-affine problems, the intrusive approach does not yield significant speed-up \cite{romor2023non}. Conversely, the non-intrusive approaches, which are purely data-driven (with no equations requirement), ensure a significant speed-up even for non-linear/non-affine problems \cite{Tezzele2022}. We focus exclusively on a non-intrusive approach, which is discussed later in \autoref{sec:Reduced Order Modelling}. In works \cite{ohlberger2013nonlinear, Taddei2015, taddei2020registration}, incorporating the mappings to obtain transformed linear approximation subspace for the construction of efficient ROMs, the construction of mapping is separated from the construction of solution coefficients and in \cite{mojgani2017arbitrary} the mappings and solution coefficients are built simultaneously. Learning mapping and solution coefficients simultaneously lead to non-linear and non-convex optimization problems, even for linear PDEs \cite{taddei2020registration}. Furthermore, applying parametric model reduction techniques to the problems incorporating mappings is standard \cite{taddei2020registration}. Conversely, when dealing with non-linear manifolds/trial subspaces, developing rapid and reliable procedures for online computation might be difficult and require several involved strategies \cite{lee2020model, romor2023non}. \\

The previous work of NNsPOD \cite{papapicco2022neural} primarily focuses on the non-linear model reduction of the hyperbolic equations. \textcolor{Black}{It solely reports the transformation of the solution manifold by exploiting deep learning framework to incorporate automatic shift detection}, for seeking optimal mapping of numerous snapshots to the reference frame and followed by POD, to derive low-dimensional transformed linear subspace. However, it is important to emphasize that it addresses only a portion of the offline stage within the ROM framework—specifically, the non-linear model reduction using NNsPOD. The complete algorithm comprising both the offline and online stages for the construction of efficient ROM is not reported. \textcolor{Black}{In this work, we construct ROM on the transformed linear approximation subspace obtained by NNsPOD and further, we employ automatic shift detection in the online phase to map/transport back the predictions in a correct physical frame. Thereby, we propose a complete NNsPOD-ROM algorithm comprising offline-online stages for the construction of accurate and efficient ROMs. } \\

The main contribution of our work is as follows:

\begin{itemize}

    \item We construct reduced order models on the low-dimensional transformed linear approximation subspace, obtained by performing non-linear model reduction using NNsPOD. This employs automatic shift detection by exploiting neural networks to derive optimal mapping of the numerous snapshots to the reference frame and aims to accelerate the KnW decay for advection-dominated problems \cite{papapicco2022neural}.  

    \item We propose the complete NNsPOD-ROM algorithm comprising both offline and online stages, where we integrate the general framework of non-intrusive parametric ROMs within the NNsPOD approach to construct accurate and efficient ROMs. We came across that the works reported in \cite{kovarnova2022shifted, burela2023parametric} propose a workflow for the development of ROM using sPOD comprising offline-online stages, in which neural networks are employed for mapping parameters to modal coefficients and also for the mapping of shifts, particularly for the predictions in the online stage. The framework proposed in \cite{Zorawski2024} is inspired by NNsPOD, it solely addresses model reduction applied to the manufactured test cases, and the complete workflow comprising offline-online stages for the construction of  ROMs is not reported. Contrary, we propose a complete NNsPOD-ROM algorithm comprising offline-online stages.
       
    \item We employ automatic shift detection in the online stage, to transport/map back the predicted field values at the reference frame (where ROM is constructed) to the physical frame/location in the physical space. 

    \item We evaluate the proposed methodology through various numerical experiments, by assessing its performance on a 1D linear advection equation, the popular higher order benchmark test case - 2D isentropic convective vortex \cite{spiegel2015survey, wang2013high} and the 2D two-phase flow test case with a focus on capturing the interface, whose physics is inspired by an industrial scenario. Both the 2D problems are governed by hyperbolic conservation laws. Further, we highlight that we focus on time-dependent problems where there is slow decay of KnW due to the advection nature of the problem.

\end{itemize}

This paper is organised as follows, in \autoref{sec:model_problem} we present the nature of the full-order model considered in this study; in \autoref{sec:Reduced Order Modelling}, we outline the general framework for the construction of non-intrusive data-driven ROMs; in \autoref{subsec:Neural Network shift augmented transformation}, we introduce the neural network shift augmented transformation approach; in \autoref{NNsPOD-ROM}, we discuss the implementation of the complete NNsPOD-ROM algorithm comprising both offline-online stages; in \autoref{sec:Numerical Experiments} we test our methodology by applying to advection-dominated problems governed by hyperbolic PDEs, by first applying to 1D linear advection equation to illustrate the complete workflow of the algorithm, then the higher order benchmark testcase - 2D isentropic convective vortex and 2D two-phase flow test cases; and finally the conclusions, perspectives and future directions are discussed in \autoref{sec:Conclusion}. 

%% file: sections/methodology.tex
\color{black}

\section{Full order model}
\label{sec:model_problem}

Let us consider the abstract problem defined by the following general time-dependent PDE for an unknown space-time dependent scalar field $\mathbf{u (x, \textit{t})}$:

\begin{equation}
    \frac{\partial \mathbf{u}}{\partial t} + \mathcal{L}(\mathbf{u}) = \mathcal{S}(\mathbf{u}, t), \quad \mathbf{x} \in \Omega, \quad t \in T
    \label{eqn:model_problem}
\end{equation}

where $\frac{\partial \mathbf{u}}{\partial t}$  represents the rate of change of the scalar field with respect to time, $\mathcal{S}(\mathbf{u}, t)$ is a source term (e.g., external forcing or a sink term), and $\mathcal{L}(\mathbf{u})$ is a differential operator that governs the evolution of $\mathbf{u}$ due to various processes such as advection, diffusion and/or other dynamical interactions. The operator $\mathcal{L}(\mathbf{u})$ governs the dynamics of the problem and can be either linear or non-linear. In many physical systems, $\mathcal{L}(\mathbf{u})$ combines both an advection term, $\mathbf{v} \cdot \nabla \mathbf{u}$ and a diffusion term, $\mathcal{D}\nabla^2 {\mathbf{u}}$:

\begin{equation}
    \mathcal{L}(\mathbf{u}) = \mathbf{v} \cdot \nabla \mathbf{u} + \mathcal{D}\nabla^2 {\mathbf{u}}.
\end{equation}

\noindent Here, $\mathbf{v}(\mathbf{x}, t)$ represents the velocity field and $\mathcal{D}$ is the diffusivity. In advection-dominated problems, the dynamics are primarily governed by the advection term $\mathbf{v} \cdot \nabla \mathbf{u}$, which describes the transport of the scalar field $\mathbf{u}$ by the velocity field $\mathbf{v}$, rather than by diffusion or other processes. In such cases, the diffusion term is negligible compared to the advection term. \\

To solve such time-dependent PDE (\ref{eqn:model_problem}) with given initial and boundary conditions, one can use any discretization method mentioned earlier (like FDM, FVM) to obtain semi-discretized equations or ordinary differential equations (ODEs) and further, the choice of time integration schemes, like Runge-Kutta schemes. However, it is important to note that solving such problems is often computationally expensive, and the solution manifold of advection-dominated problems (where diffusion is negligible) exhibits slow KnW decay. \\

Let $\mathcal{M}$ represent the solution manifold of the advection-dominated problem (\ref{eqn:model_problem}), consisting of all the solutions $\boldsymbol{\mathbf{u}}\left(\mathbf{x}, t\right)$, where $\mathbf{x}$ lies in the spatial domain $\Omega$ and \textit{t} spans the temporal domain $T$. The solution manifold is defined as:

\begin{equation}
    \mathcal{M}:=\left\{\boldsymbol{\mathbf{u}}\left(\mathbf{x}, t\right) \hspace{2mm} | \hspace{2mm} \mathbf{x} \in \Omega, t \in T \right\} 
\end{equation}

As we are dealing with time-dependent problems, the FOM solutions database is represented as a collection of time-snapshot pairs, denoted by \( \{ t_i, \mathbf{u}_i \}_{i=1}^{N_T} \), where \( t_i \in T \subset \mathbb{R}^{N_T} \) corresponds to \( N_T \) discrete time instances, and \( \mathbf{u}(t_i) \in \mathbb{R}^{N_h} \) represents the solution state with \( N_h \) degrees of freedom at each time instance $t_i$. These snapshots are organized column-wise to form the snapshot matrix \( X \in \mathbb{R}^{N_h \times N_T} \), which can be expressed as:

\begin{equation}
X=\left[\begin{array}{cccc}
\mid & \mid & & \mid \\
\mathbf{u}(t_1) & \mathbf{u}(t_2) & \ldots & \mathbf{u}(t_{N_{T}}) \\
\mid & \mid & & \mid
\end{array}\right] \in \mathbb{R}^{N_h \times N_T}.
\label{eqn:snapshots_matrix}
\end{equation}

\color{black}

\section{Reduced order modelling}\label{sec:Reduced Order Modelling}

In this section, we present the general framework for non-intrusive data-driven ROMs based on proper orthogonal decomposition (POD) \cite{Tezzele2022}, a popular linear compression method that has been extensively used in the ROM community over recent decades, especially for applications in computational fluid dynamics \cite{rozza2022advanced}. \\

\subsection{Proper orthogonal decomposition}\label{subsec:SVD}

We obtain the linear approximation subspace by performing linear dimensionality reduction using POD. The snapshots $\mathbf{{u}_\textit{i}}$ can be expressed by a linear combination of global basis functions as shown in (\ref{eqn:linearassumption}), approximating the snapshots by few global basis functions (spatial modes), $\varphi_{k}(\mathbf{x}) \in \mathbb{R}^{N_h}$ and few solution coefficients (temporal coefficients), $c_i^k(t)$. The linear approximation subspace $\mathcal{V}_{N_T}= \left\{\varphi_1,\dots,\varphi_{N_T} \right\} \in \mathbb{R}^{{N_h}\times N_T} $ is fixed, independent of the $\mathbf{{u}}_i$ to be approximated. To obtain the spatial modes $\varphi_{k}(\mathbf{x})$, we compute POD modes using singular value decomposition (SVD) \cite{PODstewart1993early} as shown in (\ref{eqn:svd}).

\begin{equation}
\mathbf{{u}}_i = \sum_{k=1}^{N_T} \varphi_k(\mathbf{x}) c_i^k(t) \approx \sum_{k=1}^r \varphi_k(\mathbf{x}) c_i^k(t) \quad \forall i \in[1, \ldots, {N_T}], \quad r \ll N_T, 
\label{eqn:linearassumption}
\end{equation}

\begin{equation}
    \underbrace{X}_{N_h \times N_T}=\underbrace{\mathbf{U}}_{N_h \times N_h} \underbrace{\boldsymbol{\Sigma}}_{N_h \times N_T} \underbrace{\mathbf{V}^T}_{N_T \times N_T} \Rightarrow 
    \underbrace{X}_{N_h \times N_T} \approx \underbrace{\boldsymbol{\hat{\mathbf{U}}}}_{N_h \times r} \underbrace{\boldsymbol{\hat{\Sigma}}}_{r \times N_T} \underbrace{\boldsymbol{\hat}{\mathbf{V}}^T}_{N_T \times N_T}, 
    \label{eqn:svd} 
\end{equation}
    
\noindent where, $\mathbf{U}$ and $\mathbf{V}$ are the two orthonormal matrices and hence, unitary matrices,  $\mathbf{U}^T \mathbf{U} = \mathbf{V}^T \mathbf{V} = \mathbb{I}$, where $\mathbb{I}$ is the identity matrix.
Columns of $\mathbf{U}$ are the left singular vectors, also called POD modes and columns of $\mathbf{V}$ are the right singular vectors of the transformed snapshot matrix $\mathcal{X}$ respectively. \color{black} $\mathbf{\Sigma}$ is the diagonal matrix with \textit{L} non-zero real singular values, arranged in descending order $\sigma_{1} \geq \sigma_{2} \geq \cdots \geq \sigma_{L}>0$, which indicates the energy contribution of the corresponding modes. As stated in the Schmidt–Eckart–Young–Mirsky theorem \cite{eckart1936approximation}, our goal is to approximate the columns of $\mathcal{X}$ by identifying a low-dimensional subspace spanned by $ r < L $ orthonormal left singular vectors: 

\begin{equation}
    \mathcal{V}_r = \left\{ \varphi_1,\dots,\varphi_{r} \right\} \in \mathbb{R}^{{N_h}\times r}. 
\end{equation}

This is achieved by selecting the first \textit{r} singular values and their corresponding left singular vectors, also known as POD modes. The resulting POD modes satisfy the following minimization problem \cite{quarteroni2015reduced}:

\begin{equation}
 \textcolor{black}{\min _{\mathcal{V}_{r}} \left\|\mathcal{X} - \mathbf{U}_r \mathbf{U}_r^T \mathcal{X}\right\|_F,}
\end{equation}

where $ \left\|  \cdot \right\|_F $ denotes a frobenius norm, $\mathcal{V}_r$ is the linear subspace of reduced dimension \textit{r}, $\mathcal{X}$ is the transformed snapshots matrix and $\mathbf{U}_r$ consists of the first \textit{r} columns of $\mathbf{U}$ (POD modes). Thus, the reduced linear subspace is the set of vectors that minimize the distance between the snapshots and their projection onto the subspace spanned by the POD modes. \\

The energy retained by the first \textit{r} POD modes is given by (\ref{eqn: POD tolerance}), where \( r \) is chosen such that the \( E(r) \)  exceeds a specified threshold. 
\begin{equation}
    E(r) = \frac{\sum_{i=1}^r \sigma_i^2}{\sum_{i=1}^{L} \sigma_i^2}.
    \label{eqn: POD tolerance}
\end{equation}

\color{black}
Once we have the $r$ most energetic POD modes (spatial modes), we compute the modal coefficients (temporal coefficients) as shown in (\ref{eqn:coeff matrix}), by projecting the snapshots matrix onto the POD modes. The initial snapshot set can only be reconstructed using these coefficients by following the procedure outlined in (\ref{eqn:linearassumption}). To predict accurate snapshots for new, unseen parameters, we need to build parameters to the modal coefficients map,  which is discussed in the next subsection. Then, just by simple matrix multiplication, we can obtain the transformed snapshots of interest. 

\begin{equation}
 c(t_i) = \mathbf{U}_r^T \hspace{1mm} \mathbf{{u}}(t_i)
\label{eqn:coeff matrix}
\end{equation}

\subsection{Solution manifold approximation}

There are several choices available for solution manifold reconstruction described by the modal coefficients corresponding to the initial database, in our case the snapshots matrix $X$. Here, the main aim is to construct a regression $s: T \mapsto \mathbb{R}^r $ which approximates the map $f$ as in (\ref{mapping_param_to_coeff}) given by a set of $N_{train}$ input-output pairs $ \left\{t_i, c_i\right\}_{i=1}^{N_{train}}$, where $t_i$ is the parameter value associated with the $i^{th}$ snapshot and $c_i$ are the corresponding temporal coefficients obtained by (\ref{eqn:coeff matrix}).

\begin{equation}
    f : t \in T \mapsto c \in \mathbb{R}^{r}. 
    \label{mapping_param_to_coeff}
\end{equation}

Then this constructed regression model is used to predict the snapshot $\tilde{\mathbf{u}}^*$ for a new unseen input parameter $t^*$ by computing, 

\begin{equation}
    {\mathbf{u}}^* = \mathbf{U}_r \hspace{1mm}s(t^*).
\end{equation}

To device this regression $s$ for the parameters to coefficients mapping, there are several techniques \cite{Tezzele2022}, namely linear interpolation, non-linear Gaussian process regression (GPR), radial basis function (RBF)  interpolation \cite{Xiao2015, taddei2020registration}, multi-fidelity methods, inverse distance weighting (IDW) and artificial neural networks (ANNs) \cite{hesthaven2018non}, to name few. In this work, we mainly use radial basis function (RBF) interpolation to build parameters to modal coefficients map.  \\

\textcolor{black}{The solution manifold of advection-dominated problems exhibits slow KnW decay, meaning it requires a high-dimensional linear approximation subspace to accurately represent the solutions. This creates a significant challenge for efficient model reduction, as the computational efficiency of ROMs depends on the reduced dimensionality of the linear approximation subspace. To circumvent this issue, we employ a manifold transformation technique that accelerates the KnW decay and enables the construction of efficient and accurate ROMs using linear methods. Further, we highlight that this transformation technique doesn't require any a-priori knowledge of the physical properties (i.e., characteristic advection velocities) and the underlying equations of the given problem.}

\section{Neural network shift augmented transformation}\label{subsec:Neural Network shift augmented transformation}

In this section, we introduce the neural network shift augmented transformation technique, which leverages deep learning architectures to employ automatic shift detection by learning a shift operator. We exploit this operator to align numerous snapshots to a reference frame, thereby accelerating the KnW decay and obtaining a low-dimensional transformed linear subspace that effectively approximates the solutions of the advection-dominated problems. Furthermore, we discuss the construction of efficient and accurate non-intrusive ROMs by following the procedure outlined in \autoref{sec:Reduced Order Modelling}, utilizing the transformed linear subspace. \\

The neural network shift augmented transformation technique introduced within the NNsPOD framework \cite{papapicco2022neural} is inspired by the shifted-POD (sPOD) approach \cite{reiss2018shifted}. In sPOD, snapshots are spatially transported to a reference frame using a shifting operator $\mathbf{\mathcal{T}_b}$, which acts on the space-time dependent field $\mathbf{u}( \mathbf{x}, t)$ as shown in (\ref{eqn: shift operator}). The resulting shifted snapshot, $\mathbf{\tilde{u}}$, corresponds to a transformed spatial distribution obtained by applying the shifting quantity that is proportional to the advection velocity, \textbf{b}.

\begin{equation}
    \mathbf{\tilde{u} = \mathcal{T}_{b} u(x, \textit{t}) = u(x - b\textit{t}, \textit{t})}, \quad \mathbf{x} \in \Omega, \quad t \in T. 
    \label{eqn: shift operator}
\end{equation}

 The application of sPOD relies on comprehensive knowledge of the advection velocities of a given problem, being the significant limitation of this approach \cite{papapicco2022neural}. Although the authors of \cite{reiss2018shifted} propose different methods to deduce the shift operator when the advection velocities are not known explicitly, these methods inevitably involve complexities and limitations, which are thoroughly discussed in section 2.4 of \cite{papapicco2022neural}. The NNsPOD framework \cite{papapicco2022neural, rozza2022advanced} overcomes this limitation by employing a semi-supervised learning paradigm for automatic shift detection. This approach automatically finds the optimal shifts for numerous snapshots to the reference frame by leveraging deep learning frameworks. The capability for automatic shift detection is particularly critical when shifts are not known a priori, as determining them in advance is often impractical or even infeasible for a wide range of problems.\\
 
In NNsPOD, the neural-network shift augmented manifold transformation is followed by linear compression POD, for non-linear reduction of hyperbolic equations. In this approach, the deep-learning framework \cite{papapicco2022neural} is utilized to \textcolor{black}{derive mapping between the original manifold $\mathcal{M}$ and the shifted manifold $\tilde\mathcal{{M}}$} as shown in (\ref{eqn: bijective mapping}). This method does not require any underlying equations (physical models or governing equations) and a-priori knowledge of the given problem's physical properties (i.e., advection velocities), making it purely data-driven. Here, we learn the shift operator $\mathcal{T}_{\text{shift}}$ using neural networks, which acts on the field $ \mathbf{u(x,\textit{t})}$ as shown in (\ref{eqn: NN shift operator}), where $\mathcal{T}_{\text{shift}}: \Omega \times T \rightarrow \Omega$ is the shifting quantity that is not a fixed value linearly depending on time as in (\ref{eqn: shift operator}), but instead it depends on space and time \cite{papapicco2022neural, rozza2022advanced}.

\begin{equation}
    \mathcal{C}^{1} \ni \mathcal{T}_{\text{shift}} :  {\mathcal{M}} \mapsto \tilde{\mathcal{M}}, 
    \label{eqn: bijective mapping}
\end{equation}

\begin{equation}
     \tilde{\mathcal{M}} := \left\{ \mathbf{u(x - \mathcal{T}_{\text{shift}}(x,\textit{t}), \textit{t})}, \quad  \mathbf{x} \in \Omega, \quad t \in T \right\}.  
    \label{eqn: NN shift operator}
\end{equation}

In this work, we employ an artificial neural network (ANN), specifically a multi-layer perceptron (MLP) referred to as ShiftNet, to learn the shift operator for the given problem. In principle, one can opt to choose different architectures depending on the complexity of the problems at hand. The learning process is built on the set of snapshots as a discrete field, as defined in (\ref{eqn:Snapshots-in-discrete fields}), where $\mathbb{V}$ is the discretized domain, and $t_i \in T $ with $N_T$ number of time instances at which snapshots are collected. \textcolor{black}{The primary objective is to transform the snapshots effectively by minimizing the Euclidean distance between them and the selected reference configuration, \( \mathbf{u}_{\text{ref}} \). The choice of the reference configuration is arbitrary, as it is selected from the set of snapshots. For this task, we define the loss function $\mathcal{L}_{\text{ShiftNet}}$, as shown in (\ref{eqn:shifted-loss-function}).}     \\

After the shift detection, in a discrete context, the transformation \( \mathbf{x} - \mathcal{T}_\text{shift}(\mathbf{x}, t) \) may not lie within the discrete space \( \mathbb{V} \). This necessitates an additional step to calculate the mean difference between the snapshots and the reference configuration. \textcolor{black}{Consequently, another crucial step in this transformation technique is reconstructing the transformed snapshots field values within the same discrete space \( \mathbb{V} \), also referred to as physical space \cite{papapicco2022neural, reiss2018shifted}. To facilitate the reconstruction of field values in the same discrete space, a neural network termed InterpNet is integrated into the algorithm. InterpNet is specifically devised to learn an operator that serves as an interpolator during the training phase of this transformation technique. The loss function $\mathcal{L}_{\text{InterpNet}}$ considered for InterpNet is the mean squared error (MSE) function, as shown in (\ref{eqn:interpNet-loss-function}), where $\mathbf{u}_{\text{ref}}$ represents the field values of the reference snapshot, while $\mathbf{x} \in $ \( \mathbb{V} \) denotes the coordinates of the cell centroids or nodes in the discretized domain $\mathbb{V}$ and $N_h$ is the number of degrees of freedom. The reference snapshot is chosen arbitrarily from the snapshot set. 
}

\begin{equation}
    \left\{\mathbb{V} \ni \mathbf{u}\left( t_i \right)\right\}_{i=1}^{N_T}= \left\{ \mathbf{u}\left(\mathbb{V}, t_i\right)\right\}_{i=1}^{N_T},
    \label{eqn:Snapshots-in-discrete fields}
\end{equation}  

\begin{equation}
    \mathcal{L}_{\text{InterpNet}}=\frac{1}{N_h} \sum_{j=1}^{N_h}\left\| \mathbf{u}_{\text{pred}}(\mathbf{x}_j) - \mathbf{u}_{\text{ref}}(\mathbf{x}_j) \right\|^2_2, \quad \mathbf{x} \in \mathbb{V}, 
    \label{eqn:interpNet-loss-function}
\end{equation}

\begin{equation}
\mathcal{L}_{\text{ShiftNet}} = \frac{1}{N_T} \sum_{i=1}^{N_T} \left\| \mathbf{u}_{\text{ref}} - \mathbf{{u}(x - \mathcal{T}_\text{shift}(x,\textit{t}_\textit{i}), \textit{t}_\textit{i})} \right\|^2_2, \quad \mathbf{x}
\in \mathbb{V}.
\label{eqn:shifted-loss-function}
\end{equation}

In Algorithm \ref{alg:one}, we report the neural network shift augmented transformation algorithm. In this transformation technique, the automatic shift detection and the reconstruction of field values workloads are divided between two separate neural networks : 
\begin{itemize}
    \item \textbf{{ShiftNet}} learns the optimal shift operator $\mathcal{T}_\text{shift}(\textbf{x},t)$ for a given problem. This operator precisely quantifies the optimal shifts, which results in shifted space that transports/maps numerous snapshots to the reference frame. To learn this shift operator we minimize the loss function $\mathcal{L}_\text{ShiftNet}$ given in (\ref{eqn:shifted-loss-function}).  
     
    \item \textbf{{InterpNet}} learns the reference configuration based on the grid distribution $\mathbf{x} \in \mathbb{V}$ (physical space). To achieve this, we define the loss function $\mathcal{L}_{\text{InterpNet}}$ shown in (\ref{eqn:interpNet-loss-function}). Once trained, InterpNet serves as an interpolator in the ShiftNet training process, to compute transformed snapshot field values in each shifted space.
\end{itemize}

 \noindent \textbf{Training neural networks:} InterpNet and ShiftNet are trained sequentially, with InterpNet being trained first, followed by ShiftNet. InterpNet learns to reconstruct the field values of the reference configuration, and its forward map is subsequently used in ShiftNet training to reconstruct field values in each shifted space. Therefore, effective optimization of the InterpNet loss function, $\mathcal{L}_{\text{InterpNet}}$, is crucial. \\

The input to InterpNet is $\mathbf{x}$, the spatial coordinates of cell centroids or nodes in the discretized domain $\mathbb{V}$ of the selected reference configuration, while the output is the corresponding field values $\mathbf{u}_\text{ref}$. As mentioned earlier, the reference configuration is chosen arbitrarily from the training set. Upon achieving the required criterion by minimizing the loss function $\mathcal{L}_\text{InterpNet}$, the operator $\mathcal{T}_\text{InterpNet}$ is obtained, which functions as an interpolator in ShiftNet training. \\

For ShiftNet, the inputs are $(\mathbf{x}, t_i)$, where $t_i$ represents the time instances at which snapshots are collected, and the output is the spatial coordinates $\mathbf{x}$ (physical space). During training, for each shifted space, the field values are computed using the $\mathcal{T}_\text{InterpNet}$ operator. The loss function $\mathcal{L}_{\text{ShiftNet}}$ is minimized to determine the optimal shift operator $\mathcal{T}_\text{Shift}$. It is important to note that a single shift operator is trained using the complete training dataset. \\

\noindent \textbf{Manifold transformation:} After training neural networks, the shift operator $\mathcal{T}_\text{Shift}$ is used to compute the optimal shifted space ${\mathbf{\tilde{x}}_{i}}^{\text{ref}}$, aligning multiple snapshots to the reference frame. The optimal shifted space is obtained by performing the operation ${\mathbf{\tilde{x}}_{i}}^{ref}$ = $\mathbf{x}$ - $\mathcal{T}_{\text{shift}}(\mathbf{x}, t_i)$, a procedure referred to as \textit{automatic shift detection}. After transporting/mapping all snapshots to the reference frame utilising the optimal shifted space ${\mathbf{\tilde{x}}_{i}}^{\text{ref}}$, the shifted/transformed snapshots field values must be reconstructed in physical space $\mathbf{x}$ (where ${\mathbf{x}} = \mathbf{x}_{\text{ref}} $), in the same grid distribution. This step is crucial in achieving sharp KnW decay and is performed using a linear interpolator $\mathcal{I}$. After this step, all the transformed snapshots are at the reference frame and in the same grid distribution $ \mathbf{x} \in \mathbb{V}$.

\RestyleAlgo{ruled}
\SetKwComment{Comment}{\(\triangleright\) }{ }

\begin{algorithm}[H]
\DontPrintSemicolon
\KwData{Inputs ($\mathbf{x}, t_i$) and corresponding snapshots $\mathbf{u}_{i} \in \textbf{X}$, where i = 0,...,$N_{\text{train}}$; reference configuration index ; InterpNet and ShiftNet architecture and respective criterion $\epsilon_{\text{interp}}, \epsilon_{\text{shift}}$ for loss functions. }
\KwResult{Transformed snapshots at reference frame and $\mathcal{T}_{\text{shift}}$ optimal shift operator. }
    
    \While{$ \mathcal{L}_{\text{InterpNet} } > \epsilon_{\text{interp}}$}{
        \textbf{InterpNet}.forward ; \\
        compute  $ \mathcal{L}_{\text{InterpNet}} $ ; \\
        \textbf{InterpNet}.backward ; \\
    }
    
    $\mathcal{T}_{\text{InterpNet}}$ = $ \textbf{InterpNet}$.forward; \Comment*[r]{reconstructs field values}
    
    \For{$\mathbf{u}(t_i)\in \mathbf{X}$, where i = 0 $\mathbf{to}$ $N_{\text{train}}$}{
        \While{$\mathcal{L}_{ShiftNet} > \epsilon_{\text{shift}}$}{
            ${\mathbf{\tilde{x}}_{i}}$ = $\mathbf{x}_{i}$ - \textbf{ShiftNet}.forward ; \Comment*[r]{x is physical space (x = x$_{\text{ref}}$)}
            $ \tilde{\mathbf{u}}_i=\mathcal{T}_{\text{interp}}({\mathbf{\tilde{x}}_{i}}); $ \Comment*[r]{reconstruct field values at each shifted space} 
            compute  $\mathcal{L}_{ShiftNet}$ ; \\
            \textbf{ShiftNet}.backward ; \\
        }
        
        $\mathcal{T}_{\text{shift}}$=\textbf{ShiftNet}.forward; \Comment*[r]{optimal shift operator}
        $\tilde{\mathbf{x}}_{i}^{ref}$ = $\mathbf{x}$ - $\mathcal{T}_{\text{shift}}$ ; \Comment*[r]{shifted space that transport snapshot to reference frame}
        $\mathcal{L} = \mathcal{I}(\tilde{\mathbf{x}}_{i}^{ref},\tilde{\mathbf{u}}_{i})$ ; \Comment*[r]{interpolator}
        $\tilde{\mathbf{u}}_{i}^{ref} = \mathcal{L}(\mathbf{x}) $; \Comment*[r]{transformed snapshots in physical space} 
    }
\caption{Neural network shift augmented transformation}
\label{alg:one}
\end{algorithm}

\section{NNsPOD-ROM algorithm}\label{NNsPOD-ROM}

In this section, we present the implementation of the complete NNsPOD-ROM algorithm, which includes both the offline and online stages for constructing efficient and accurate ROMs for advection-dominated problems. As outlined in Algorithm \ref{alg:two}, the NNsPOD-ROM algorithm is divided into two stages: 1) FOM Pre-processing (Offline Stage): This stage involves the neural network shift-augmented transformation (which includes the training of neural networks and manifold transformation) and the construction of the ROM on the transformed linear subspace, and 2) FOM Post-processing (Online Stage): This stage focuses on the testing and prediction phase. \\

\noindent \textbf{FOM Pre-processing (offline stage):}
At first, we collect input parameters $(\mathbf{x}, t_i)$ and their corresponding snapshots $\mathbf{u}_i$, which are required for training. We then select a reference configuration, specify the neural network architectures for InterpNet and ShiftNet, and define their stopping criteria. Initially, we obtain the transformed snapshots matrix $\mathcal{X}$ and the shift operator $\mathcal{T}_\text{shift}$ by following Algorithm \ref{alg:one}, which reports neural network shift augmented transformation. After this step, all the transformed snapshots are at the reference frame and in the same grid distribution $ \mathbf{x} \in \mathbb{V}$. Finally, we perform linear dimensionality reduction on the transformed snapshots matrix, to achieve sharp KnW decay and obtain a low-dimensional transformed linear approximation subspace. We obtain the \textit{r} number of most energetic POD modes to construct ROM and follow the framework outlined in \autoref{sec:Reduced Order Modelling} to construct non-intrusive or data-driven ROMs. \\

\noindent \textbf{FOM Post-processing (online Stage):} Initially, the predicted snapshot field values from the constructed ROM are in the reference frame (where the ROM is built). Therefore, shifting or transporting these field values back to the correct frame in physical space is necessary. This is accomplished by computing the shifted space ${\mathbf{\tilde{x}}_{i}^{shift}} = \mathbf{x} + \mathcal{T}_{\text{shift}}(\mathbf{x}, t_i)$, which transports/maps the predicted snapshot field values to the correct frame. Here, we highlight that we are employing \textit{automatic shift detection} for predictions in the online phase. The predicted snapshot field values are transported to the correct frame utilising the shifted space ${\mathbf{\tilde{x}}_{i}^{shift}}$. Subsequently, the transported field values are then reconstructed in physical space $\mathbf{x}$ using a linear interpolator $\mathcal{I}$. The predicted snapshot field values intensity variations in the correct frame/location in the physical space $\mathbf{x} \in \mathbb{V}$. \\

\RestyleAlgo{ruled}
\SetKwComment{Comment}{\(\triangleright\) }{ }

\begin{algorithm}[H]
\DontPrintSemicolon
\KwData{
  Inputs (\(\mathbf{x}, t_i\)) and corresponding snapshots \(\mathbf{u}_{i} \in \textbf{X}\), where \(i = 0, \dots, N_{\text{train}}\); \\
  \(\textit{r}\) number of POD modes; reference configuration index; \\
  InterpNet and ShiftNet architecture and respective criteria \(\epsilon_{\text{interp}}, \epsilon_{\text{shift}}\) for loss functions.
}
\KwResult{NNsPOD-ROM prediction of snapshots for given parameter \(t_i\)}

\underline{FOM pre-processing: offline stage} \\
\vspace{1mm}
Obtain transformed snapshots matrix \(\mathcal{X}\) and shift operator \(\mathcal{T}_{\text{shift}}\) using Algorithm \ref{alg:one}; \\
\vspace{1mm}
POD = POD(rank = r); \Comment*[r]{returns \(r\) POD modes}
s = POD.singular\_values; \Comment*[r]{get singular values}
\(\tilde{\mathbf{u}}^{{ref}}_{i}\) = ROM(\(\mathbf{\mathcal{X}}\), POD, \(\bullet\)); \Comment*[r]{construct ROM, \(\bullet\) = Linear, GPR, RBF, or ANN}

\underline{FOM post-processing: online stage} \\   
\vspace{1mm}
\For{\(t_{i}\), where \(i = 0\) \(\mathbf{to}\) \(N_{\text{test}}\)}{
    \(\tilde{\mathbf{u}}^{ref}_{i} = \text{ROM.predict}(t_i)\); \Comment*[r]{predictions in reference frame}
    \({\mathbf{\tilde{x}}_{i}^{shift}} = \mathbf{x} + \mathcal{T}_{\text{shift}}\); \Comment*[r]{shifted space transports back to physical frame}
    \(\mathcal{L} = \mathcal{I}({\mathbf{\tilde{x}}_{i}^{shift}}, \tilde{\mathbf{u}}^{ref}_{i})\); \Comment*[r]{interpolator}
    \({\mathbf{u}}_i^* = \mathcal{L}({\mathbf{x}})\); \Comment*[r]{predicted field values in physical frame}
}
Get \({\mathbf{u}}_i^*\) in physical frame and compute relative \(L_2\)-norm error; \\
\({\mathbf{u}}_{\text{new}} = \text{ROM.predict}(t_{\text{new}})\); \Comment*[r]{predictions for unseen new parameters}

\caption{NNsPOD-ROM}
\label{alg:two}
\end{algorithm}

\begin{figure}[ht!]
    \centering
    \includegraphics[scale=0.45]{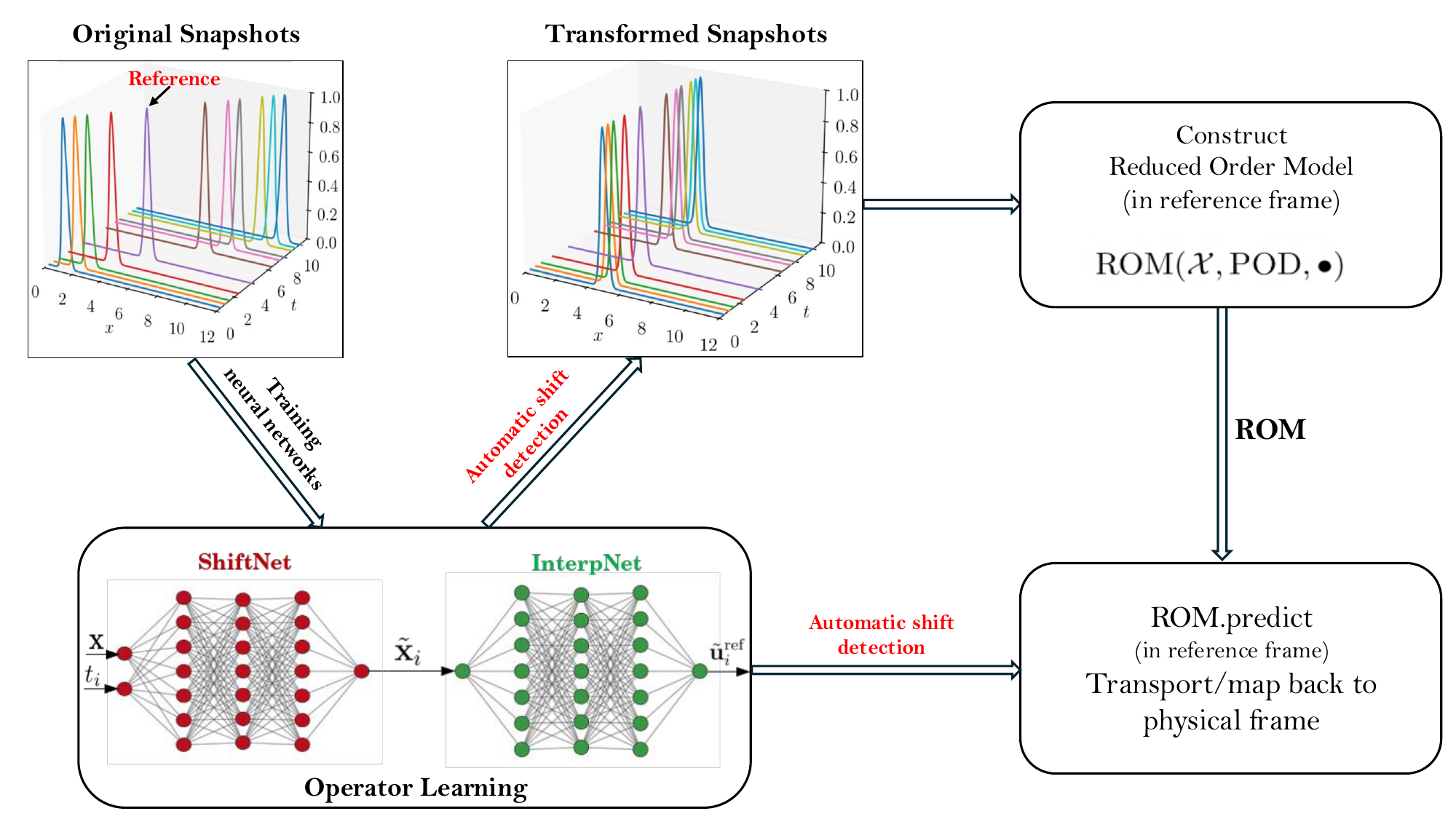}
    \caption{Schematic representation of the complete workflow for constructing NNsPOD-ROMs. The offline stage consists of training neural networks, manifold transformation and construction of non-intrusive ROM. In the online stage, automatic shift detection is employed to transport back the predicted snapshot field values to the correct frame in physical space. To illustrate the workflow, we consider 1D travelling waves, referred to as the Original Snapshots. During the training phase, we first minimizing the $\mathcal{L}_\text{InterpNet}$ (\ref{eqn:interpNet-loss-function}) of InterpNet to learn the $\mathcal{T}_\text{Interp}$ operator, which functions as interpolator to reconstruct the reference configuration field values. Subsequently, we train ShiftNet, where the $\mathcal{T}_\text{Interp}$ operator is used to reconstruct field values in the shifted space (new grid distribution). The loss function $\mathcal{L}_\text{ShiftNet}$, as defined in Equation \ref{eqn:shifted-loss-function}, is minimized to learn the optimal shift operator $\mathcal{T}_\text{ShiftNet}$. Further to perform manifold transformation, the optimal shift operator $\mathcal{T}_\text{ShiftNet}$ (automatic shift detection) is employed to transport/map numerous snapshots to the reference frame. Then the transformed snapshots' field values are reconstructed in discrete physical space $\mathbf{x}$ (referred to as Transformed snapshots) using a linear interpolator. This process accelerates KnW decay and results in a low-dimensional transformed linear approximation subspace. The non-intrusive ROM is built on this low-dimensional transformed linear subspace. In the online stage, the predicted snapshot field values are initially at the reference frame (where ROM is constructed), and therefore, automatic shift detection is employed to transport/map back the field values to the physical frame. The transported field values are reconstructed in physical space using a linear interpolator. At this stage, the predicted snapshot field values intensity variations are in the physical frame.}
    \label{fig:Schematic_NNsPOD_ROM}
\end{figure} 

To assess the accuracy of the predictions, we compute the relative $L_2$ error of the predicted snapshots $\mathbf{u}^*$, defined as in (\ref{eqn: l2 error}). Here, $\mathbf{u}_{i}$ represents the original or true field value, and ${\mathbf{u}}_i^*$ is the ROM-predicted field value at the $i^{\text{th}}$ cell.

\begin{equation}
\mathcal{L}_{2}= \frac{\left\| \mathbf{u} -  {{\mathbf{u}}}_i^* \right\|}{\left\| \mathbf{u}  \right\|} =  \frac{\sqrt{\sum_{i=0}^{N_h}\left(\mathbf{u}_{i} - {{\mathbf{u}}}_i^* \right)^2}}{\sqrt{\sum_{i=0}^{N_h}\left(\mathbf{u}_{i}^2\right)}} . 
\label{eqn: l2 error}
\end{equation}

The implementation of both Algorithm \ref{alg:one} and Algorithm \ref{alg:two} is carried out using the \texttt{EZyRB} package \cite{demo18ezyrb}, a Python library for non-intrusive data-driven model order reduction of parameterized problems and here, the deep neural networks are implemented using the \texttt{PyTorch} package \cite{paszke2017automatic}. A schematic representation of the complete workflow of the NNsPOD-based ROM algorithm is shown in \autoref{fig:Schematic_NNsPOD_ROM}. \\

%% file: sections/num_exp.tex
\section{Numerical Experiments}\label{sec:Numerical Experiments}

In this section, we evaluate the effectiveness of the proposed NNsPOD-ROM methodology by applying it to the construction of efficient and accurate ROMs for advection-dominated problems. The methodology is demonstrated by considering three test cases: a 1D linear advection equation, a popular higher-order methods benchmark case—a 2D Isentropic convective vortex and a 2D Two-phase flow, focusing on capturing the interface.

\color{black}{\subsection{1D Linear advection equation}

We consider the 1D linear advection equation (\ref{eqn:linearTransportEquation}), where \( u(x, t) \) represents the scalar quantity (i.e, concentration, temperature, etc.) being transported as a function of space \( x \in \Omega \) and time \( t \in T \). The constant \( c \) is the advection velocity representing the rate at which the quantity is transported in the \( x \)-direction. 

\begin{equation}
\frac{\partial u(x, t)}{\partial t} + c \frac{\partial u(x, t)}{\partial x} = 0, \quad x \in \Omega.
\label{eqn:linearTransportEquation}
\end{equation}

For a given initial condition \( u(x, 0) = f(x) \) at \( t = 0 \), the solution to the linear advection equation is \( u(x, t) = f(x - ct) \). A commonly used initial condition in transport problems is Gaussian function (\ref{eqn:gaussian_function}) \cite{fresca2021comprehensive, Geelen2023}, parameterized by its mean \( \mu \) and standard deviation \( \sigma \).  In this study, we utilize the linear advection equation to describe the evolution of a Gaussian wave. The initial condition is defined as:

\begin{equation}
u(x, 0) = f(x) = \exp\left(-\frac{(x - \mu)^2}{2 \sigma^2}\right),
\label{eqn:gaussian_function}
\end{equation}

where mean \( \mu \) represents the initial position of the wave, while standard deviation \( \sigma \) determines the width of the wave. In this case, both values are kept constant, \( \mu = 0.75 \) and  \( \sigma = 2 \). The solution to the linear advection equation with this initial condition is obtained in closed form. As we are dealing with time-dependent problems, we parametrize time within the range \( t \in [0, 10]\). At each time instance \(t \), the quantity \( u(x, t) \) simply shifts along the \( x \)-axis at a constant rate \( c = 1 \) without changing its shape. We considered this test case to illustrate the complete workflow of the proposed algorithm. We generate a database consisting of 21 snapshots at uniformly sampled time instances \( t \in [0, 10]\), over 256 equidistant nodes  $ x \in [0, 12] $, as shown in \autoref{fig:1Doriginal_gaussian_snaps}. We divide 50$\%$ of the database as a training set (\autoref{fig:train database}, Left) and the rest as a testing set. \\
\color{black}

\begin{figure}[ht]
    \includegraphics[scale=0.8]{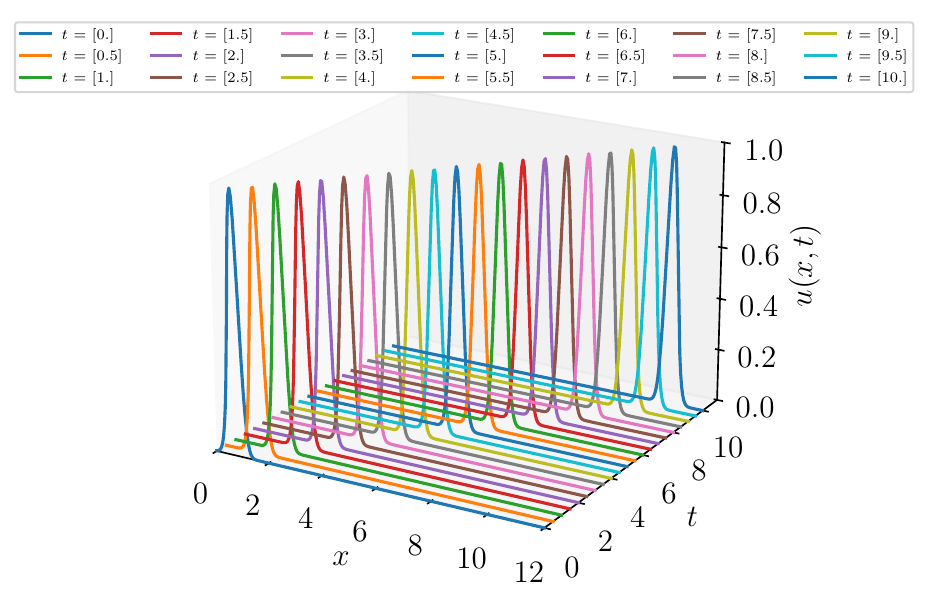}
    \centering
    \caption{Snapshots of the 1D linear advection equation collected at 21 uniformly sampled time instances, \( {t} \in [ 0, 10 ]\), over 256 equidistant spatial nodes $ x \in [0, 12] $. }
    \label{fig:1Doriginal_gaussian_snaps}
\end{figure}

For training neural networks, the snapshot corresponding to the reference parameter $\text{t}_{\text{ref}}$ = 4 from the training set, is selected as the reference configuration. The settings for the two neural network architectures used in this test case are shown in \autoref{table:1Dwave}. InterpNet learns to reconstruct the field values of the reference snapshot as shown in \autoref{fig:InterNetLearning and ShiftNetLearning} (left). This is utilised in the ShiftNet training to learn the optimal shift operator, which yields the shifted space that transports/maps the original snapshot field values to reference frame/configuration. In \autoref{fig:InterNetLearning and ShiftNetLearning} (right), we show the original snapshot field values mapped to the reference frame utilizing the respective shifted space for the parameter $t$ = 10, is labelled as Shifted, the original snapshot of the same parameter $t$ = 10, is labelled as physical and the reference configuration is labelled as reference. It can be noticed in \autoref{fig:InterNetLearning and ShiftNetLearning} (right) that, the shifted snapshot field values don't lie in the same discretized domain $x \in [0, 12]$ or physical space. As mentioned in the previous section, the shifted/transformed snapshot field values must be reconstructed in physical space using interpolation to achieve sharp KnW decay. In \autoref{fig:train database}, the comparison between the original snapshots of the training set (left) and the reconstructed shifted snapshots at the reference frame in physical space are shown. \\

\begin{figure}[h!]
    \begin{subfigure}[h]{0.5\textwidth}
         \includegraphics[width=\textwidth]{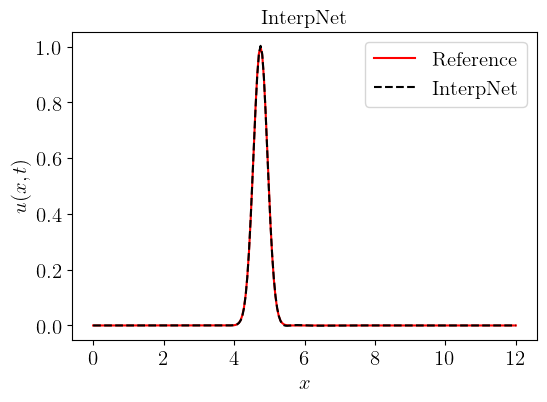}
    \end{subfigure}
    \begin{subfigure}[ht]{0.5\textwidth}
        \includegraphics[width=\textwidth]{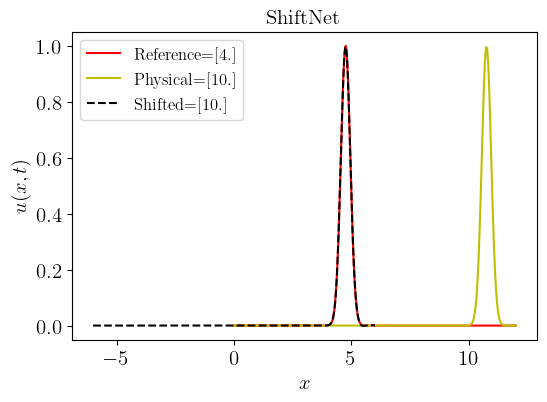}
    \end{subfigure}
\caption{For 1D linear advection equation, \textbf{Left}: Showing trained InterpNet, able to reconstruct the field values of reference snapshot corresponding to time instance, $t_\text{ref} = 4$ and \textbf{Right}: The original snapshot at time t = 10 (labelled as physical), the same transported to the reference frame utilising the shifted space (labelled as shifted), using optimal shift obtained by shift operator $\mathcal{T}_\text{shift}$ and the reference frame is also shown for comparison.}
\label{fig:InterNetLearning and ShiftNetLearning}
\end{figure}

\begin{figure}[h!]
    \includegraphics[scale=0.7]{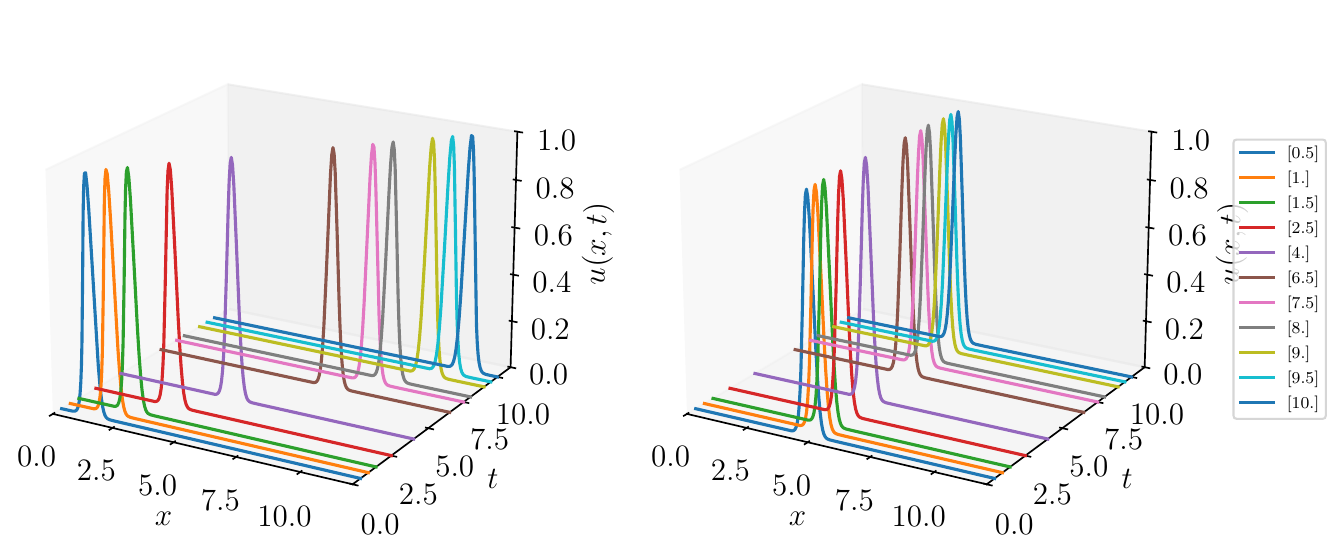}
    \centering        
    \caption{For 1D linear advection equation, \textbf{Left}: showing training set snapshots in physical space $x$ and \textbf{Right}: the shifted snapshots at the reference frame in physical space after neural network shift augmented transformation.}
\label{fig:train database}
\end{figure}

\begin{figure}[h!]
    \begin{subfigure}[h!]{0.5\textwidth}
         \centering
         \includegraphics[width=\textwidth]{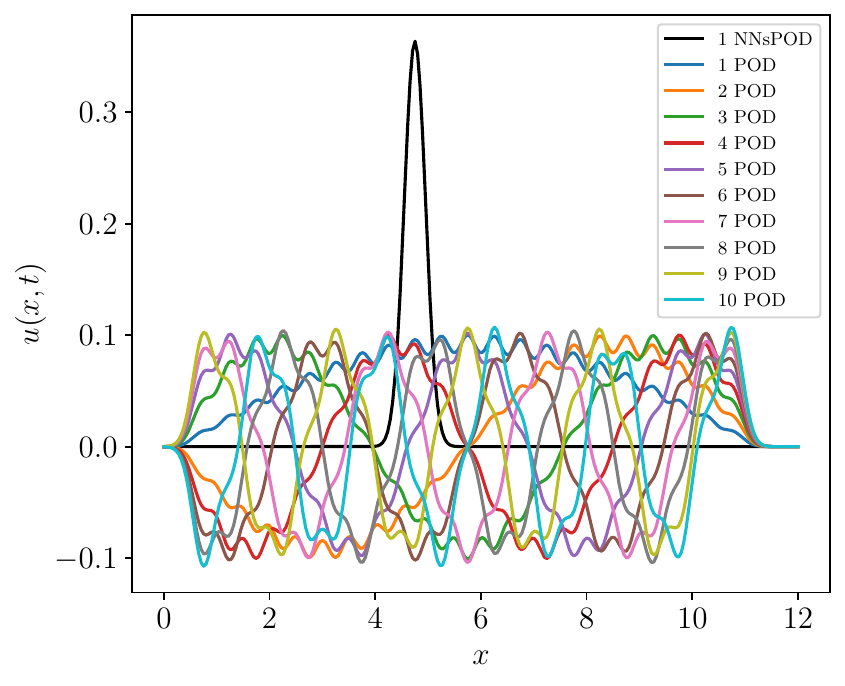}
         \label{fig:1DPODmode}
    \end{subfigure}
    \begin{subfigure}[h!]{0.5\textwidth}
        \centering
        \includegraphics[width=\textwidth]{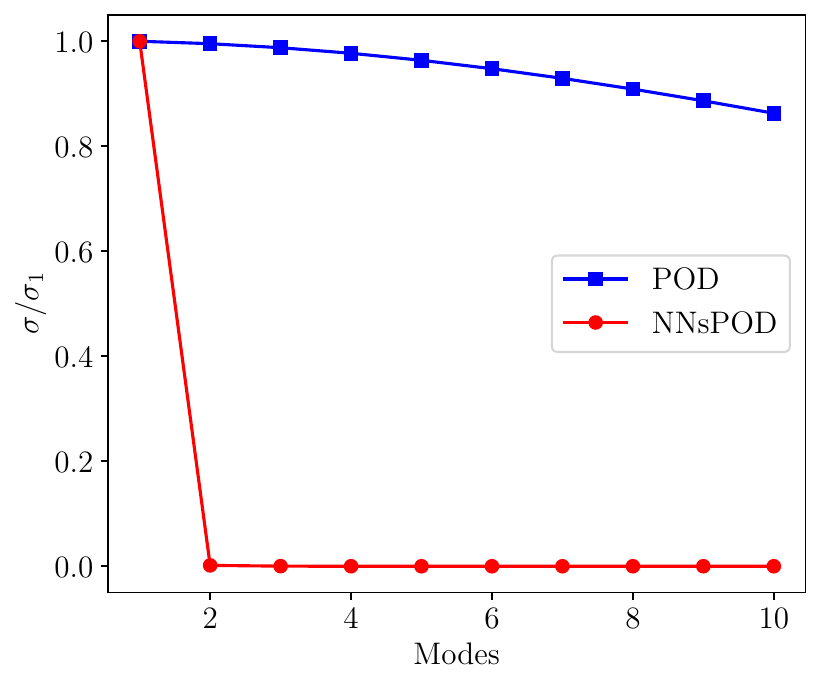}
        \label{fig:1DShiftedSingularValues}
    \end{subfigure}
\caption{For 1D linear advection equation, \textbf{Left}: Comparison of the first NNsPOD mode, obtained via neural-network shifted augmented manifold transformation, with all POD modes of the full-order model solutions. and \textbf{Right}: Singular value decay comparison, showing POD (without transformations) and NNsPOD (with the neural-network shifted augmented manifold transformation)}
\label{fig:pod and singular values}
\end{figure}

After the neural-network shift augmented transformations of the training set database, all the snapshots are transported/mapped to the reference frame, as shown in \autoref{fig:train database}. The low-dimensional linear approximation subspace is then obtained by performing singular value decomposition (SVD) of the transformed snapshot matrix $\mathcal{X}$. \autoref{fig:pod and singular values} demonstrates the effectiveness of NNsPOD over transitional POD, showing the sharp singular values decay of transformed matrix $\mathcal{X}$ and the most energetic first NNsPOD mode, capturing $99.99\%$ of the cumulative energy. In contrast, achieving the same energy requires all the POD modes. Furthermore, the POD modes fail to capture the shape of the Gaussian wave, whereas NNsPOD captures the Gaussian wave accurately. To construct efficient ROMs, the modal coefficients are computed by projecting the transformed snapshot matrix $\mathcal{X}$ onto the NNsPOD subspace. To obtain fast and accurate predictions for given unseen parameters, we build parameters to modal coefficients map using radial basis function (RBF) interpolation. Then, we obtain the snapshot for unseen parameters by simple matrix multiplication. Initially, the predicted snapshots are in the reference frame. As discussed in the previous section, we employ automatic shift detection to transport/map back the predicted snapshot field values to the physical frame. \\

\color{black}
In \autoref{fig:1DPrediction_test}, the predicted snapshot field values (solid line) for test set parameters in the physical frame are compared with full-order linear advection equation solutions (dotted line). The reconstruction errors for the training and test set parameters are shown in \autoref{fig:1D_prediction_error}, with the error magnitudes exhibiting similar variations. This demonstrates that NNsPOD-RBF-ROM predicts reliable accurate snapshots for the new unseen parameter. In \autoref{fig:1D_prediction_POD_NNsPOD_comparision}, we compare the mean relative $L^2$ reconstruction error of standard POD-RBF-ROM and NNsPOD-RBF-ROM against the modes rank for test set parameters of the linear advection equation. By considering the first NNsPOD mode with a cumulative energy of $99.99\%$, the NNsPOD-RBF-ROM provides accurate predictions, outperforming standard POD-RBF-ROM. Even when all POD modes are considered, the performance of the standard POD-RBF-ROM remains significantly inadequate, being ill-suited for advection-dominated problems. \\

\color{black}
\begin{figure}[h!]
    \includegraphics[scale=0.8]{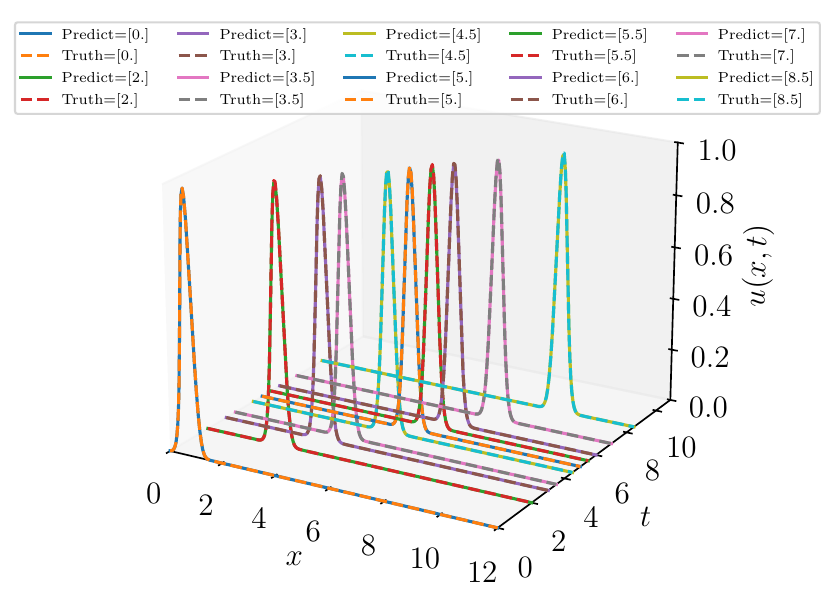}
    \centering
        \caption{Prediction of NNsPOD-RBF-ROM (considered first NNsPOD mode with 99.99$\%$ cumulative energy) for test parameters, labelled as Predict (solid lines), compared with the full-order solutions of 1D linear advection equation, labelled as Truth (dashed lines).}
    \label{fig:1DPrediction_test}
\end{figure}

\begin{figure}[h!] 
    \centering
    \includegraphics[scale=0.5]{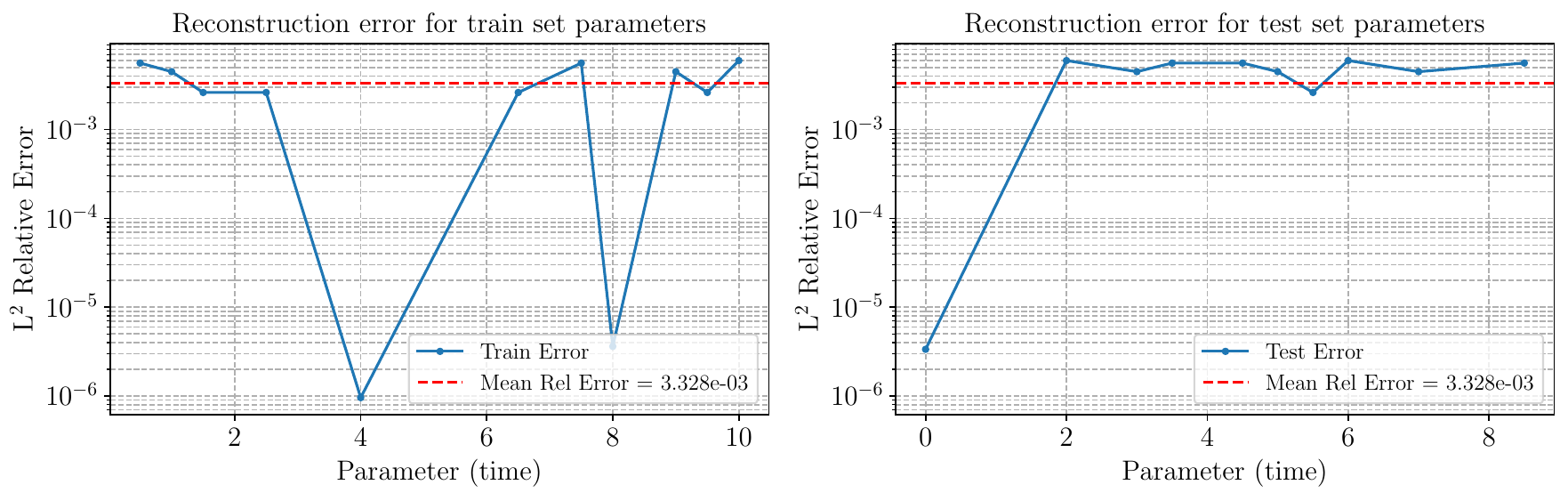}
        \caption{Relative $L^2$ reconstruction error of predictions obtained by NNsPOD-RBF-ROM for the 1D linear advection equation, shown for training set parameters (left) and test set parameters (right). The red line on the left indicates the mean training error, while the red line on the right represents the mean test error, demonstrating that the resulting ROM exhibits excellent generalization capabilities.}
    \label{fig:1D_prediction_error}
\end{figure}

\begin{figure}[h!] 
    \centering
    \includegraphics[scale=0.6]{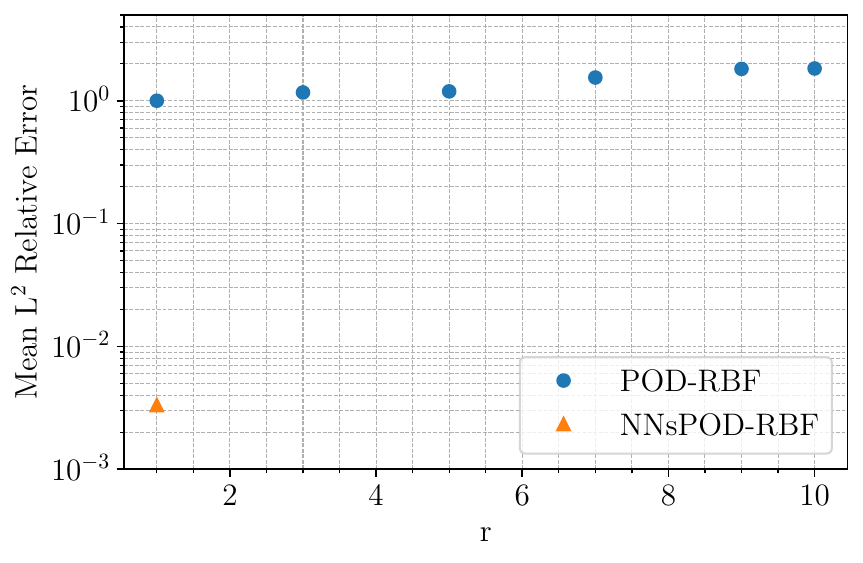}
        \caption{Mean $L^2$ relative reconstruction error comparison of POD-RBF-ROM and NNsPOD-RBF-ROM predictions for test set parameters of 1D linear advection equation.}
    \label{fig:1D_prediction_POD_NNsPOD_comparision}
\end{figure}

\subsection{2D Isentropic convective vortex}
We consider the popular higher order methods benchmark case \cite{spiegel2015survey, wang2013high}, isentropic convective vortex to test the proposed NNsPOD-ROM methodology. The governing equations considered are the 2D Unsteady Euler equations (\ref{eqn:Euler equation}), to describe the motion of the inviscid, compressible flow of our benchmark case. \\ 

\begin{equation}
\frac{\partial}{\partial t}\left[\begin{array}{c}
\rho \\
\rho u \\
\rho v \\
\ E
\end{array}\right]+\frac{\partial}{\partial x}\left[\begin{array}{c}
\rho u \\
\rho u^2 + p \\
\rho u v \\
u(E+p)
\end{array}\right]+\frac{\partial}{\partial y}\left[\begin{array}{c}
\rho v \\
\rho u v \\
\rho v^2 + p\\
v(E+p)
\end{array}\right]=0.
\label{eqn:Euler equation}
\end{equation}

where $\rho$ is density, ($u,v$) are velocity x and y component, $p$ is pressure and $E$ is total energy. The thermodynamic closure is obtained by the total energy equation (\ref{eqn: Energy equation}) and in our case we consider constant specific heat ratio, $\gamma = 1.4$ and gas constant, $R_{gas} = 287.15$ J/Kg.K.

\begin{equation}
    E = \frac{p}{\gamma-1} + \frac{1}{2}\rho ( u ^{2} + v^{2}).
    \label{eqn: Energy equation}
\end{equation}

The initial condition of the freestream flow is given by $\rho_\infty = 1, u_\infty = 0.1, v_\infty = 0, p_\infty = 1$ and the vortex is imposed in the domain by considering the analytical solution as in (\ref{Eqn: initial condition vortex}), where vortex strength \(b=0.5\), the radius of the vortex, \(r = \left[(x-x_c)^2 + (y-y_c)^2 \right]^{1/2}\) = 0.5 and vortex centre \(\left(x_c,y_c\right) = \left(5,10\right)\). There are different versions of the vortex problem available in the community \cite{spiegel2015survey, ChiWang1997ENO-WENO}, we consider the version implemented in the package HyPar 1.0 - Finite Difference Hyperbolic-Parabolic PDE Solver on Cartesian Grids \cite{Hypar1} and the results of full-order solutions are obtained using this package.  
\begin{equation}
\left\{\begin{array}{l}
    \rho = \left[ 1 - \frac{\left(\gamma-1\right)b^2}{8\gamma\pi^2} e^{1-r^2} \right]^{\frac{1}{\gamma-1}}, \\
    p = \rho^\gamma, \\
    u = u_\infty - \frac{b}{2\pi} e^{\frac{1}{2}\left(1-r^2\right)} \left(y-y_c\right), \\
    v = v_\infty + \frac{b}{2\pi} e^{\frac{1}{2}\left(1-r^2\right)} \left(x-x_c\right). \\
\end{array}\right.
\label{Eqn: initial condition vortex}
\end{equation}

We consider a fifth-order hybrid-compact weighted essentially non-oscillatory (WENO) scheme to approximate the spatial derivatives of the governing equation, WENO scheme belongs to the class of higher-order accurate numerical methods for hyperbolic PDEs and convection-dominated problems \cite{zhang2016ENO-WENO}. We consider the third-order strong-stability preserving runge-kutta (SSPRK3) method \cite{gottlieb2001SSPRK}, for time integration. \\ 

We consider a domain $\Omega = [0, 40] \times [0, 20]$ discretized into 28800 Cartesian grids, the time interval considered for the simulation is $t \in [0, 62.5]$ with a time step, $\Delta t$ = 0.00625. We collect snapshots every 100th time-step,  yielding ${N}_T = 100$ number of snapshots and leading to the solution manifold definition \( \mathcal{M}:=\left\{\boldsymbol{\rho}\left(t_i\right) \in \Omega, i=1, \ldots, 100\right\} \), consisting of density field snapshots for 100 instances. \\

\begin{figure}[h]
    \includegraphics[scale=0.5]{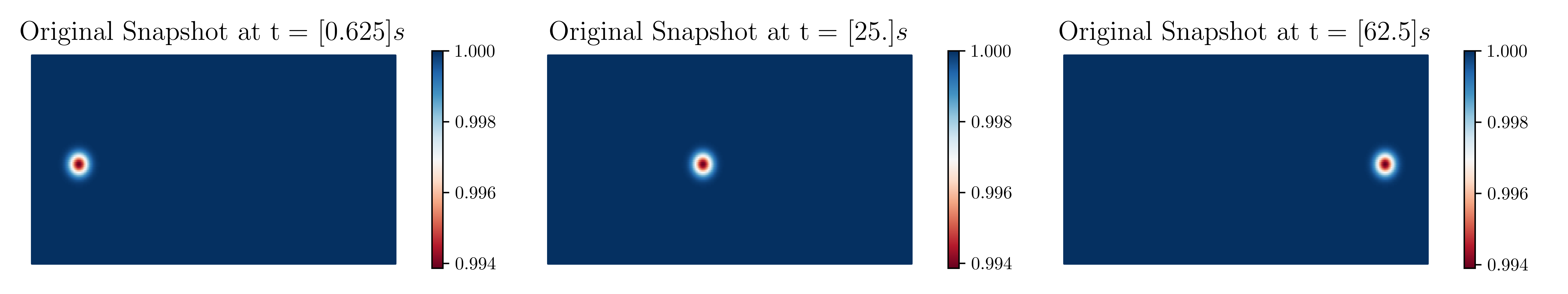}
    \centering
    \caption{FOM solution of isentropic convective vortex at time, t $\in \{0.625, 25, 62.5\}$ and for FOM pre-processing, the original snapshot at t = 0.625 (left) is considered as reference configuration.} 
    \label{fig:2Dvortex_originalSnapshots}
\end{figure}

\begin{figure}[h]
    \includegraphics[scale=0.5]{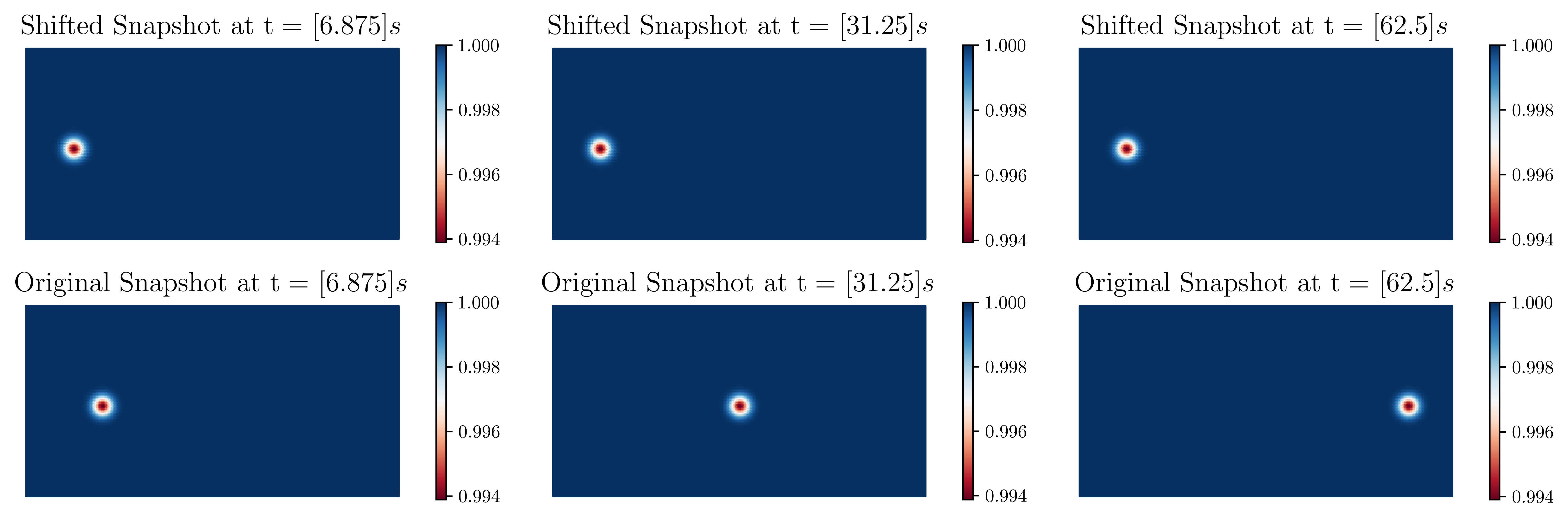}
    \centering
        \caption{Shifted snapshots (top) of isentropic convective vortex after neural-network shift augmented transformations and the respective original snapshot (bottom) at time, t $\in \{0.625, 25, 62.5\}$ from the training set.} 
    \label{fig:2Dvortex_Shifted_Snapshots}
\end{figure}

\begin{figure}[h!]
    \centering
    \begin{subfigure}[h]{0.40\textwidth}
         \centering
         \includegraphics[width=\textwidth]{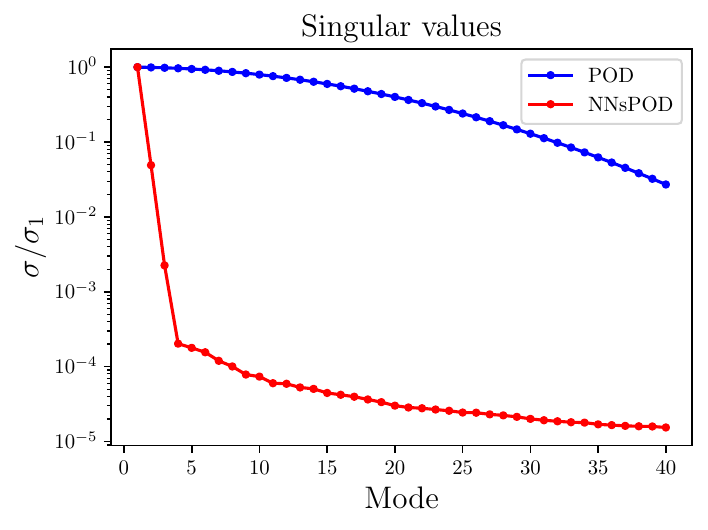}
    \end{subfigure}
    \begin{subfigure}[ht]{0.55\textwidth}
        \centering
        \includegraphics[width=\textwidth]{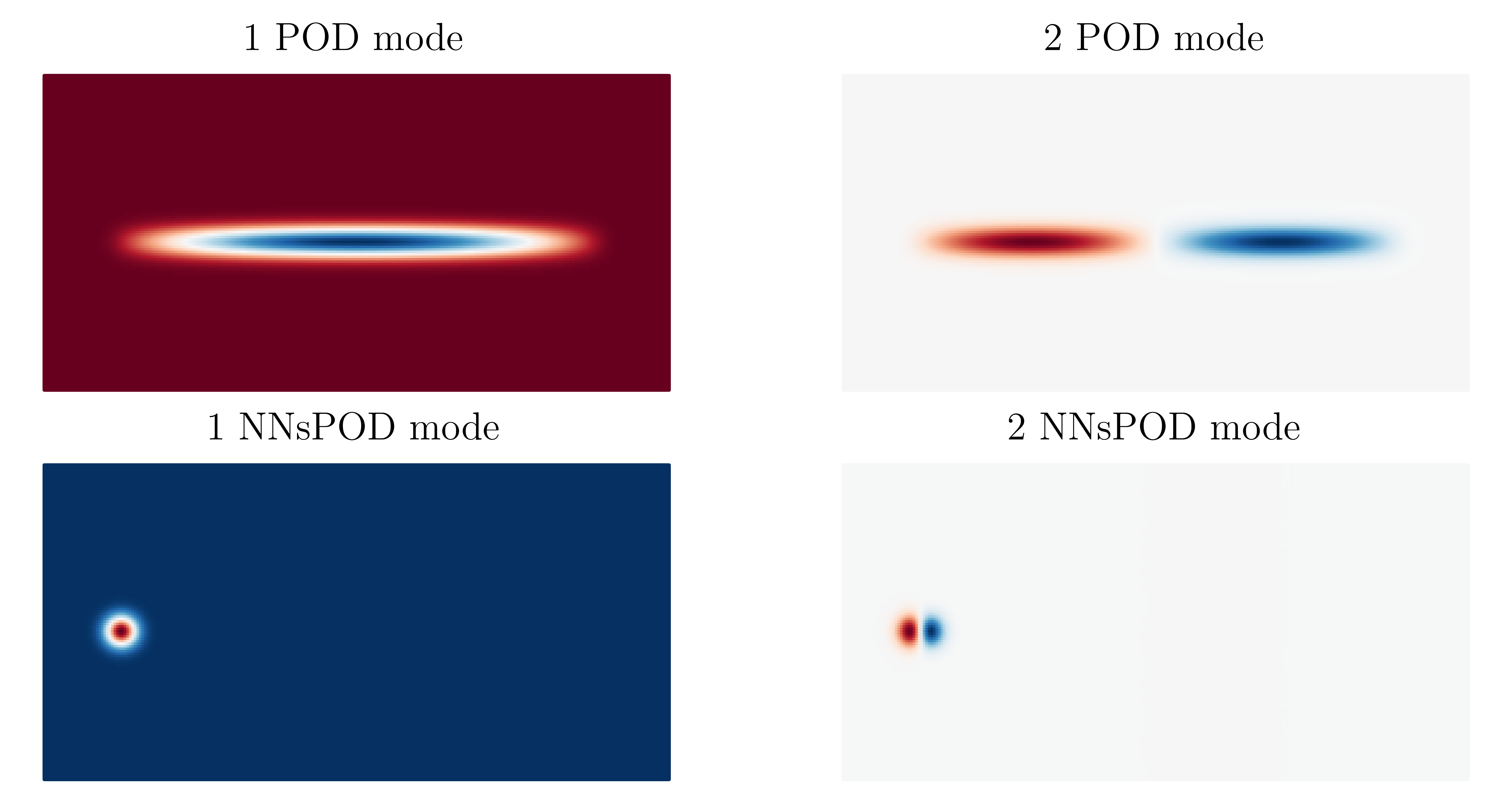}
    \end{subfigure}
    \caption{ \textbf{Left}: Comparison of singular values decay of isentropic convective vortex without the transformation, is labelled as POD and after the neural-network shifted augmented manifold transformation, is labelled as NNsPOD, and \textbf{Right}: Showing first two modes of FOM manifold without the transformations (first row - POD modes) and after the neural-network shift augmented manifold transformation (second row - NNsPOD modes).}
    \label{fig:vortex_sing_compare_POD_Mode}
\end{figure}

For this numerical experiment, we generate a database consisting of 100 snapshots for each time instance, from which we randomly divide 80$\%$ of the database into the training set and the remaining 20$\%$ into the test set. Before training, we scale the database to range between 0 and 1 using the Min-Max Scaler, reported in \autoref{subsec:scaling}. For training neural networks, we select the reference configuration parameter, t = 0.625, the corresponding snapshot is shown in \autoref{fig:2Dvortex_originalSnapshots} (left). In this case, only the x-ordinates are shifted, as the vortex is convected in the x-direction. The settings considered for two neural network architectures for this test case are as given in \autoref{table:Isentropic Convective Vortex}. InterpNet learns to reconstruct the reference snapshot field values and is utilised in ShiftNet training to learn the optimal shift operator that yields optimal shifted space, which maps the original snapshot field values to the reference frame. After the neural network shift augmented manifold transformation, we obtain shifted snapshots that are transported/mapped to the reference frame as shown in the \autoref{fig:2Dvortex_Shifted_Snapshots}. \\

We then perform linear dimensionality reduction on this transformed manifold, to obtain sharp singular values decay as shown in \autoref{fig:vortex_sing_compare_POD_Mode} (left), compared with singular values obtained without the manifold transformation. In \autoref{fig:vortex_sing_compare_POD_Mode}(right), the first two modes of the FOM manifold without transformations and after the neural network shift augmented transformations are shown. It can be observed that NNsPOD accurately captures the vortex structure, whereas the standard POD fails to do so. Furthermore, the first 4 NNsPOD modes account for $99.99 \%$ of the cumulative energy, whereas 40 POD modes are required to constitute the same energy. We consider the first 4 NNsPOD modes to build efficient NNsPOD-RBF-ROM, and its predictions (rescaled to the original range) for unseen test set parameters are shown in \autoref{fig:Vortex_prediction_test}, with the original snapshots (Truth) and the corresponding absolute error. \\

\begin{figure}[h!]
    \includegraphics[scale=0.5]{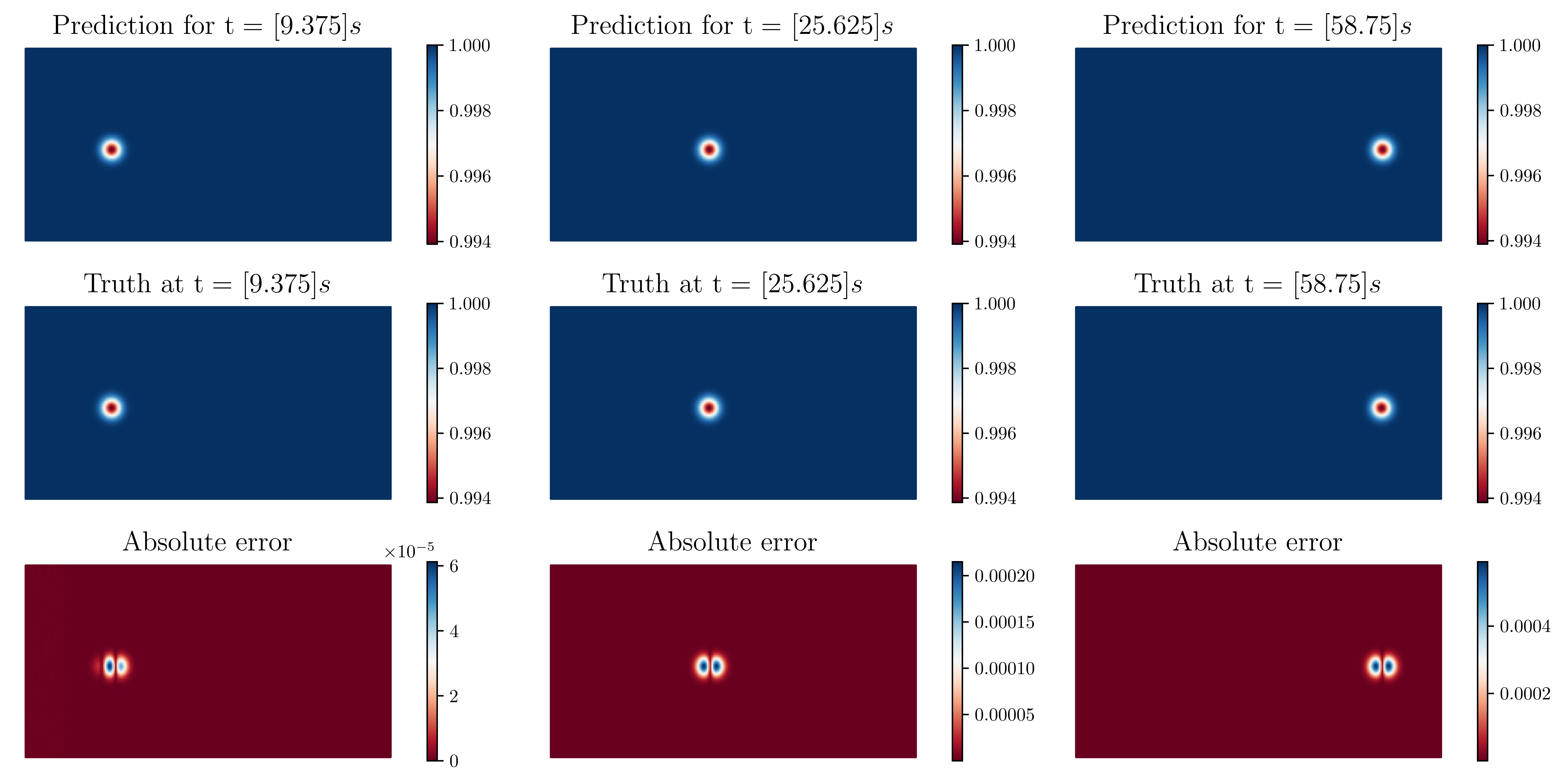}
    \centering
        \caption{Predictions of snapshots by NNsPOD-RBF-ROM of isentropic convective vortex test set parameters, showing for time instances t = $ \{9.375, 25.625, 58.75\}$s. }
    \label{fig:Vortex_prediction_test}
\end{figure}

\begin{figure}[h!] \centering
    \includegraphics[scale=0.5]{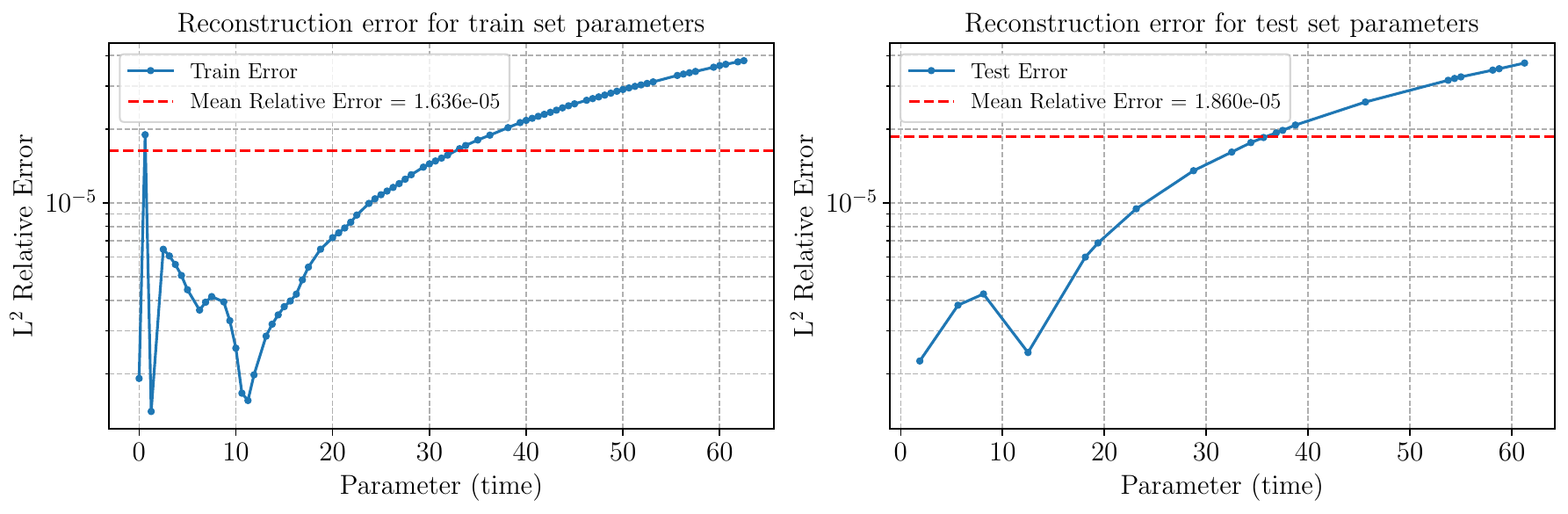}
    \centering
        \caption{Relative $L^2$ reconstruction error of predictions by NNsPOD-RBF-ROM for Isentropic convective vortex, for both training set parameters (left) and test set parameters (right).}
    \label{fig:Vortex_prediction_test_error}
\end{figure}

The reconstruction error for the training and test set parameters is shown in \autoref{fig:Vortex_prediction_test_error}. We plot the relative $L^2$ error of the predicted snapshots for each respective parameter and the mean relative $L^2$ error. The error varies in the same order of magnitude for both train and test set parameters predictions, with the error trend increasing for the parameters further from the chosen reference configuration parameter. \color{black}
In \autoref{fig:2Dvortex_prediction_POD_NNsPOD_comparision}, we compare the mean relative $L^2$ reconstruction error of the standard POD-RBF-ROM and NNsPOD-RBF-ROM to modes rank, for test set parameters of 2D isentropic convective vortex. The NNsPOD-RBF-ROM obtained by considering the first 4 NNsPOD modes that account for $99.99\%$ of the cumulative energy, provides accurate predictions and effectively captures the vortex structure. The standard POD-RBF-ROM fails to adequately represent the vortex structure, even when 10 POD modes are considered, as shown in \autoref{fig:Vortex_prediction_test_POD_RBF}. Overall, NNsPOD-RBF-ROM outperforms standard POD-RBF-ROM qualitatively and quantitatively, providing exceptional reconstruction accuracy and a more precise representation of the vortex structure via a low-dimensional transformed linear approximation subspace. \\

\begin{figure}[h!]
    \includegraphics[scale=0.5]{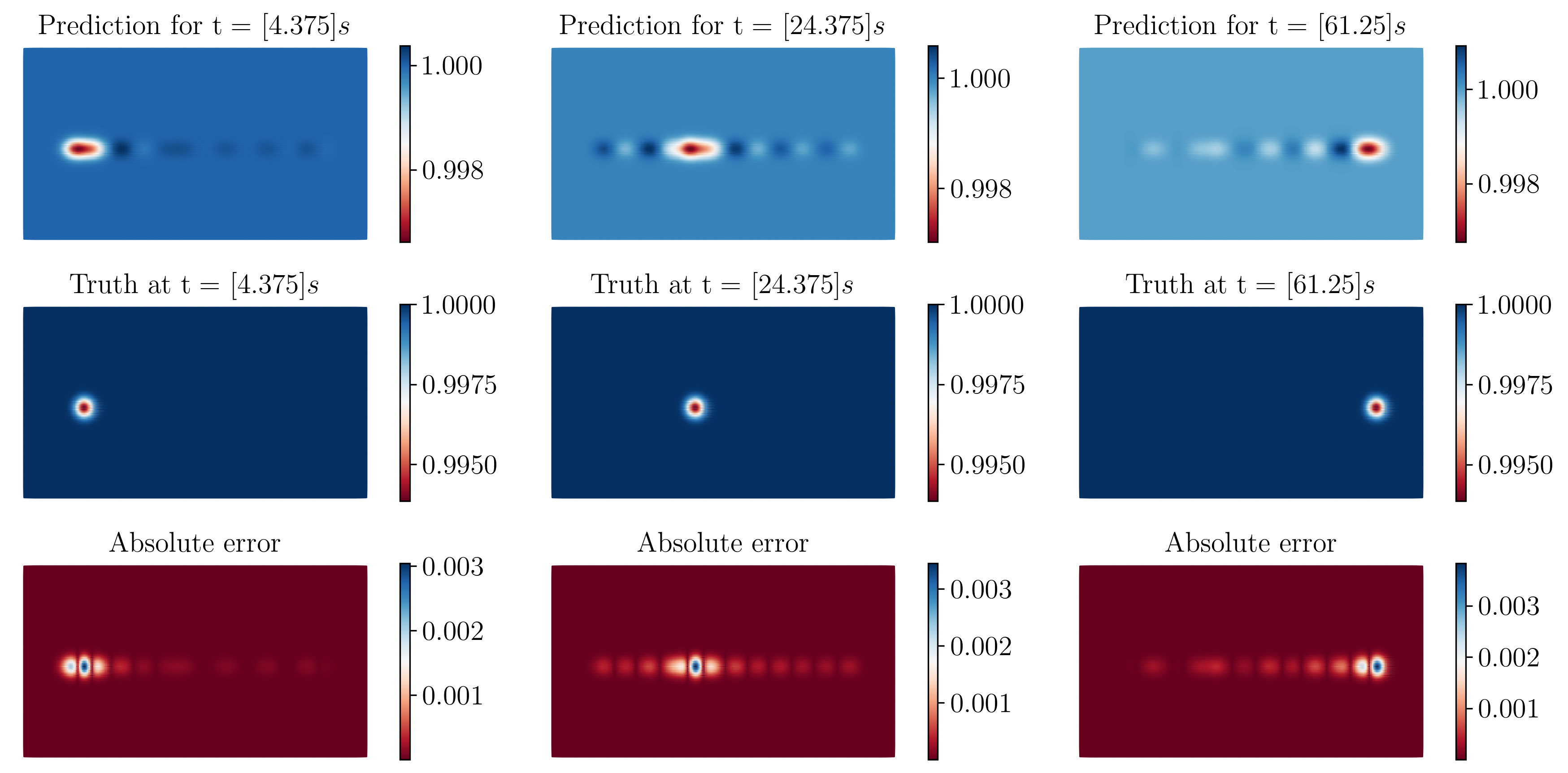}
    \centering
        \caption{Predicted snapshots by the POD-RBF-ROM for test set parameters of the isentropic convective vortex, considering 10 POD modes, shown for time instances t = $ \{4.375, 24.375, 61.75\}$ s.} 
    \label{fig:Vortex_prediction_test_POD_RBF}
\end{figure}

\begin{figure}[h!] \centering
    \includegraphics[scale=0.6]{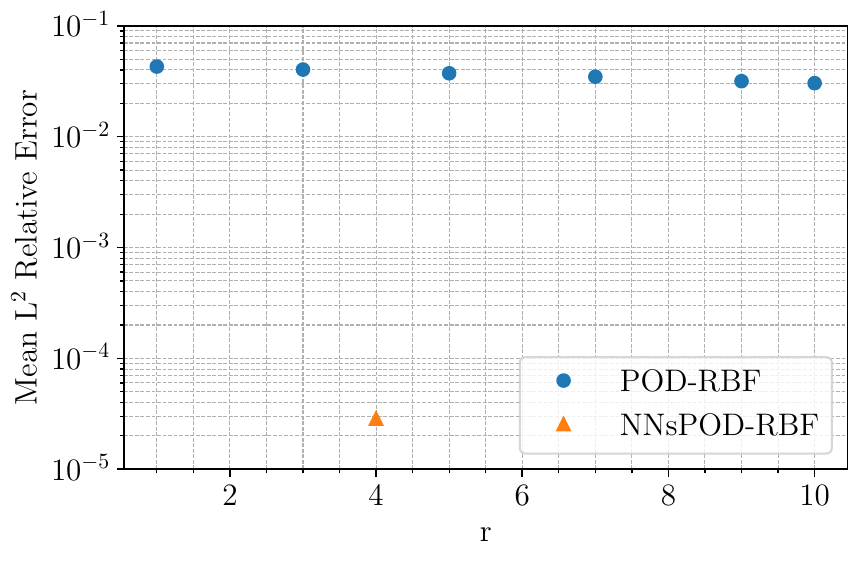}
    \centering
    \caption{Mean $L^2$ relative reconstruction error comparison of POD-RBF-ROM and NNsPOD-RBF-ROM predictions for test set parameters of 2D isentropic convective vortex.}
    \label{fig:2Dvortex_prediction_POD_NNsPOD_comparision}
\end{figure}

\color{black}
\subsection{2D Two-phase flow}

Two-phase flows are governed by the unsteady Navier-Stokes equation (\ref{Eqn:NSeqn}), with continuity and momentum equation for incompressible fluids, where $\mathbf{U}$ is velocity field, $\rho$ is density, $\textit{p}$ is pressure and $\mathbf{F_\textit{ST}}$ is the volumetric surface tension. 

    \begin{equation}
        \left\{
        \begin{array}{l}
            \partial_t(\rho) + \nabla \cdot (\mathbf{\rho U}) = 0, \\[1mm]
            \partial_t(\rho \mathbf{U})+\nabla \cdot(\rho \mathbf{U} \otimes \mathbf{U})= -\nabla p+ \nabla \cdot [ \mu (\nabla \mathbf{U} + \nabla^{T} \mathbf{U}) ]+ {\mathbf{F}_\textit{ST}}.
        \end{array}
        \right.
        \label{Eqn:NSeqn}
    \end{equation}
\vspace{3mm}
For capturing the interface between two-phase flows we consider the one-fluid models approach, which are volume-of-fluids (VOF), level-set (LS) and phase-field (PF) methods. In the one-fluid model approach, two immiscible fluids are considered as one effective fluid, in which $\mathbf{U}$ acts as the shared velocity of the two fluids $\mathbf{U} = \mathbf{U}_{1} = \mathbf{U}_{2}$. The physical properties of the fluids are averaged over the domain, which leads to constituent relations of $\rho$ and $\mu$ weighted values given by: 

\begin{equation}
\begin{array}{r}
\rho=\alpha_1 \rho_1+\left(1-\alpha_1\right) \rho_{2}, \\
\mu=\alpha_1 \mu_1+\left(1-\alpha_1\right) \mu_{2}.
\end{array}
\end{equation}

Here, we use the VOF method, in which the modified volume fraction advection equation (\ref{eqn:adv eqn}) is solved simultaneously with continuity and momentum equations. Where $\alpha \in [0,1]$ is phase fraction which takes values $\alpha = 0$ for gas, $\alpha = 1$ for liquid and at the interface, it ranges between $ 0 < \alpha < 1 $, and $\mathbf{U}_r$ is the relative velocity shown in (\ref{eqn: relative velocity}), which is perpendicular to the interface. To solve this advection equation, we use a semi-implicit MULES algorithm.  
 
\begin{equation}
\partial_t \alpha +\nabla \cdot\left(\alpha \mathbf{U} \right)+\nabla \cdot\left[\mathbf{U}_r \alpha\left(1-\alpha\right)\right]=0, 
\label{eqn:adv eqn}
\end{equation}

\begin{equation}
\mathbf{U}_r = C_\alpha|\mathbf{U}| \frac{\nabla \alpha}{\left|\nabla \alpha \right|}, \hspace{8mm} \text{with} \hspace{2mm}\mathcal{C}_\alpha \in [0,1].
\label{eqn: relative velocity}
\end{equation}

To model the turbulence, we consider $k-\epsilon$ turbulence model based on \cite{launder1983numerical}, which is a two-equation linear eddy viscosity turbulence closure model consisting of $k$ - turbulent energy and $\epsilon$ - turbulent kinetic energy dissipation rate equations: 

\begin{equation}
\left\{\begin{array}{l}
\partial_t(\rho k) = \nabla \cdot \left( \rho D_k \nabla k \right) + P - \rho \epsilon, \\ 
\partial_t(\rho \epsilon)=\nabla \cdot\left(\rho D_\epsilon \nabla \epsilon\right)+\frac{C_1 \epsilon}{k}\left(P+C_3 \frac{2}{3} k \nabla \cdot \mathbf{U}\right)-C_2 \rho \frac{\epsilon^2}{k}.
\end{array}\right.
\label{Eqn:k-epsilon}
\end{equation}

\begin{equation}
    \nu_t = C_\mu \frac{k^2}{\epsilon}.
    \label{Eqn:turbulent viscosity}
\end{equation}

\noindent where, $D_k$ - effective diffusivity for $k$, $P$ - turbulent kinetic energy production rate, $D_{\epsilon}$ - effective diffusivity for $\epsilon$, $C_1, C_2, C_3$ are model coefficients. After obtaining $k$, $\epsilon$,  values, the turbulent viscosity, $\nu_{t}$ is given by (\ref{Eqn:turbulent viscosity}), here $C_\mu$ is the model coefficient for turbulent viscosity. We consider the standard model coefficient values shown in \autoref{tab:model coefficents}. 

\begin{table}[h]
    \centering
    \begin{tabular}{ccccc}
        \toprule
        $C_\mu$ & $C_1$  & $C_2$  & $C_3$  & $\sigma$ \\
        \midrule
        0.09 & 1.44  & 1.92 & 0 & 1  \\
        \bottomrule
    \end{tabular}
    \caption{Model coefficients of $k-\epsilon$ turbulence model.}
    \label{tab:model coefficents}
\end{table}

We consider water ($\rho_{1} = 10^{3}, \nu_{1} = 10^{-6}$) and air ($\rho_{2} = 1, \nu_{2} = 1.48 \times 10^{-5}$) as two fluids for this two-phase flow test case, and consider constant inflow of water into the domain with inlet velocity $\mathbf{U}$(x, 0, 0) = (0.5, 0, 0). We consider a domain $\Omega = [-0.25, 2] \times [-0.25, 1]$ discretized into 2856 hexahedra cells, the time interval considered for the simulation is $t \in [0, 40]$ with a time step $\Delta t$ = 0.05, yielding ${N}_s = 800$ number of snapshots and leading to the solution manifold definition: 

\[
    \mathcal{M}_h:=\left\{\boldsymbol{\alpha}_h\left(t_i\right) \in \mathcal{V}_h, i=1, \ldots, 800\right\}. 
\]

We construct a database consisting of 800 snapshots, each corresponding to a specific time instance, where time is parameterized. The discretized domain used for the two-phase flow test case is shown in Figure~\ref{fig:mesh}. The inlet and outlet are highlighted in green and red, respectively. We emphasize that fluid enters only through the inlet, while the outflow at the outlet is restricted. The full-order solutions of the two-phase flow case are obtained using the interFoam solver \cite{Deshpande2012} of OpenFOAM \cite{Weller1998}, which is based on the finite volume method (FVM). In Figure \ref{fig:velocity}, we show the velocity magnitude contours at times $t \in $  \{ 10, 20, 40\}. These plots illustrate the flow complexity and help characterize the convective field. The original two-phase flow snapshots at a few time instances are shown in \autoref{fig:two_phase_originalSnapshots} and here, the inlet region at the bottom left and outlet region at the bottom right of the domain can be depicted. We are interested in capturing the interface of two-phase flow and hence, we only consider the region of a domain where the interface evolves. For this reason, we neglect inlet and outlet regions, as depicted in \autoref{fig:two_phase_Shifted_snapshots}). Then, we randomly divide 80$\%$ of the database as a training set and the rest 20$\%$ as a test set. For training neural networks, we select the last snapshot (completely filled configuration) at t = 40s as the reference configuration and the settings considered for the two neural network architectures for this test case are as reported in \autoref{table:Two-phase}. InterpNet learns to reconstruct the reference configuration, shown in \autoref{fig:two_phase_Shifted_snapshots} (bottom right) at $t=40$. In this case, since the interface primarily evolves in the y-direction, we incorporate shifts only in the y-coordinates of the spatial domain. ShiftNet learns the shift operator, which finds the optimal shift for a given parameter and is used to transport/map the original snapshot to the reference frame as shown in \autoref{fig:two_phase_Shifted_snapshots}. \\

\begin{figure}[ht!]
    \centering
    \includegraphics[width=0.6\linewidth]{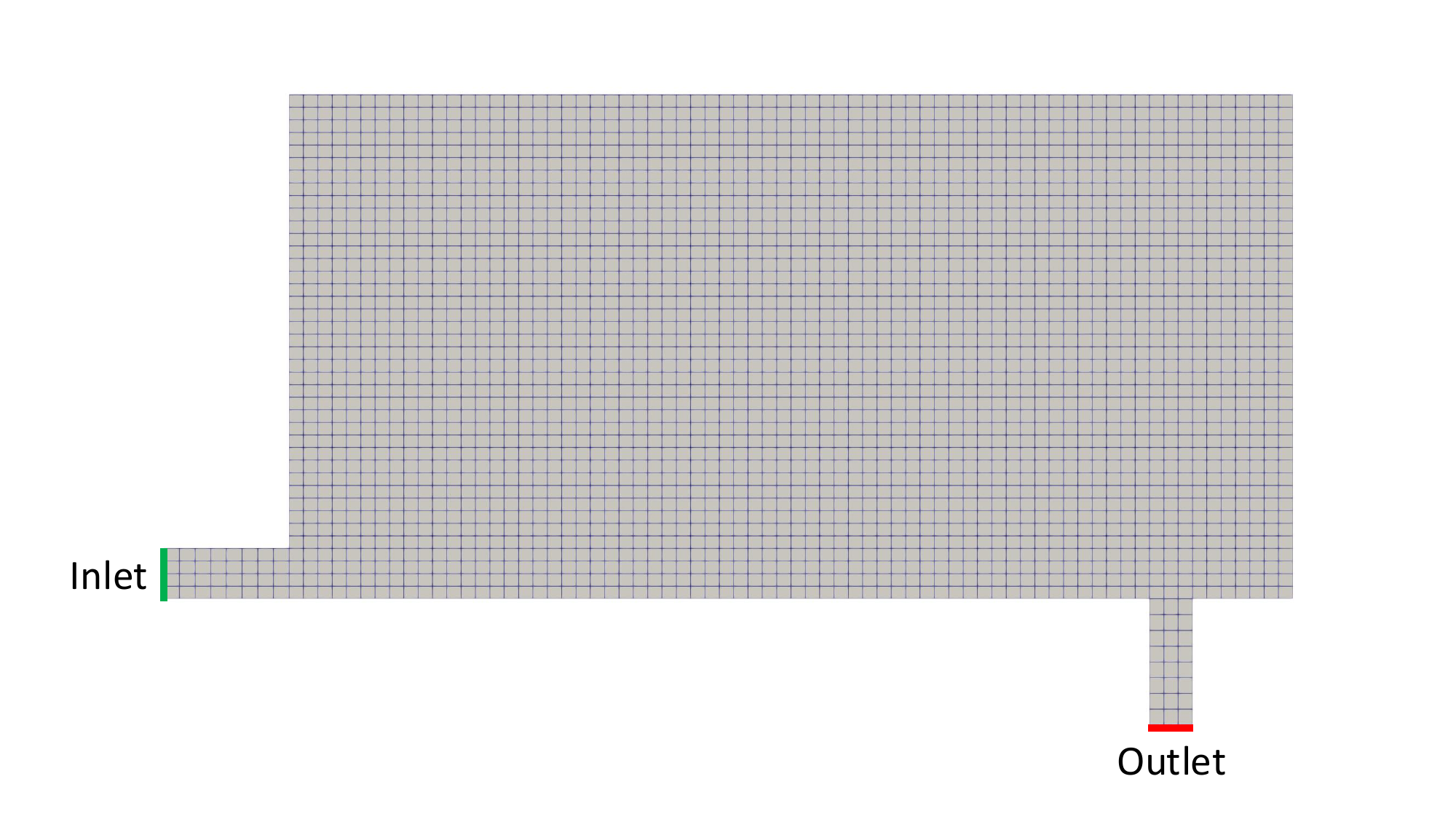}
    \caption{Discretized domain used for the two-phase flow test case. The inlet and outlet boundaries are highlighted in green and red, respectively.}
    \label{fig:mesh}
\end{figure}

\begin{figure}[ht!]
    \centering
    \includegraphics[width=1\linewidth]{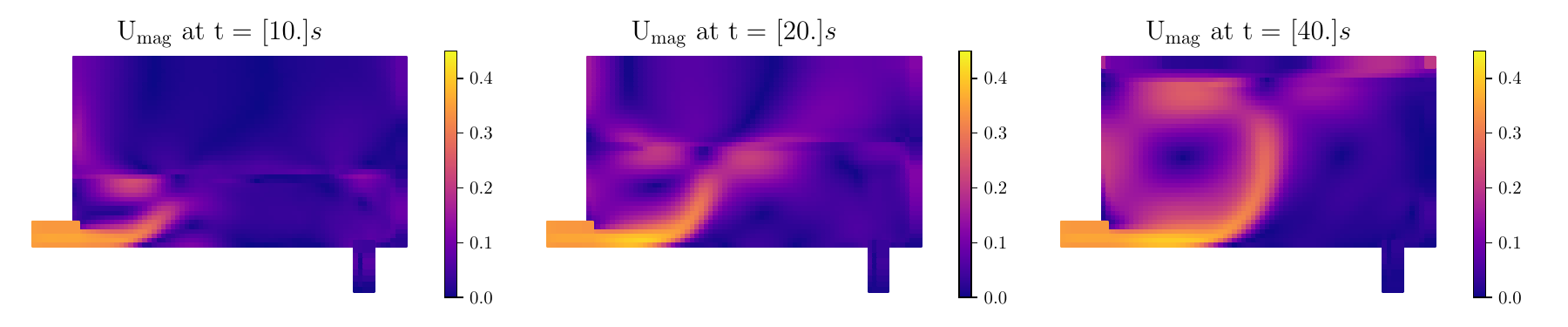}
    \caption{Velocity magnitude contour of FOM solutions for the two-phase flow at times t $\in$ \{10, 20, 40\}. 
    We highlight that fluid enters only through the inlet, while outflow through the outlet is restricted.}
    \label{fig:velocity}
\end{figure}

\begin{figure}[ht!]
    \includegraphics[scale=0.5]{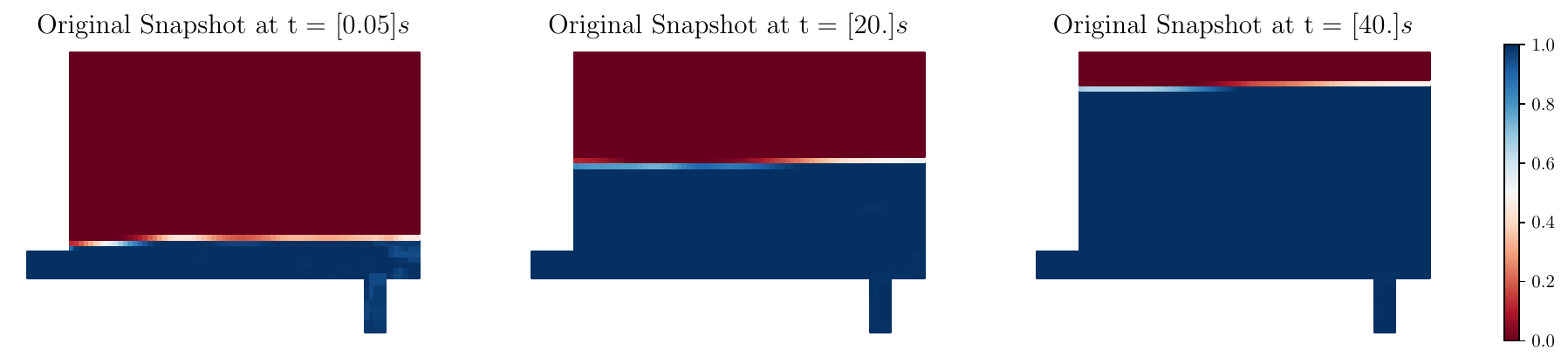}
    \centering
    \caption{FOM solution of two-phase flow at the time, t $\in \{0.05, 20, 40\}$ and for pre-processing of the solution manifold, the original snapshot at t = 40 (right figure) is considered as reference configuration.} 
    \label{fig:two_phase_originalSnapshots}
\end{figure}

\begin{figure}[ht!]
    \includegraphics[scale=0.5]{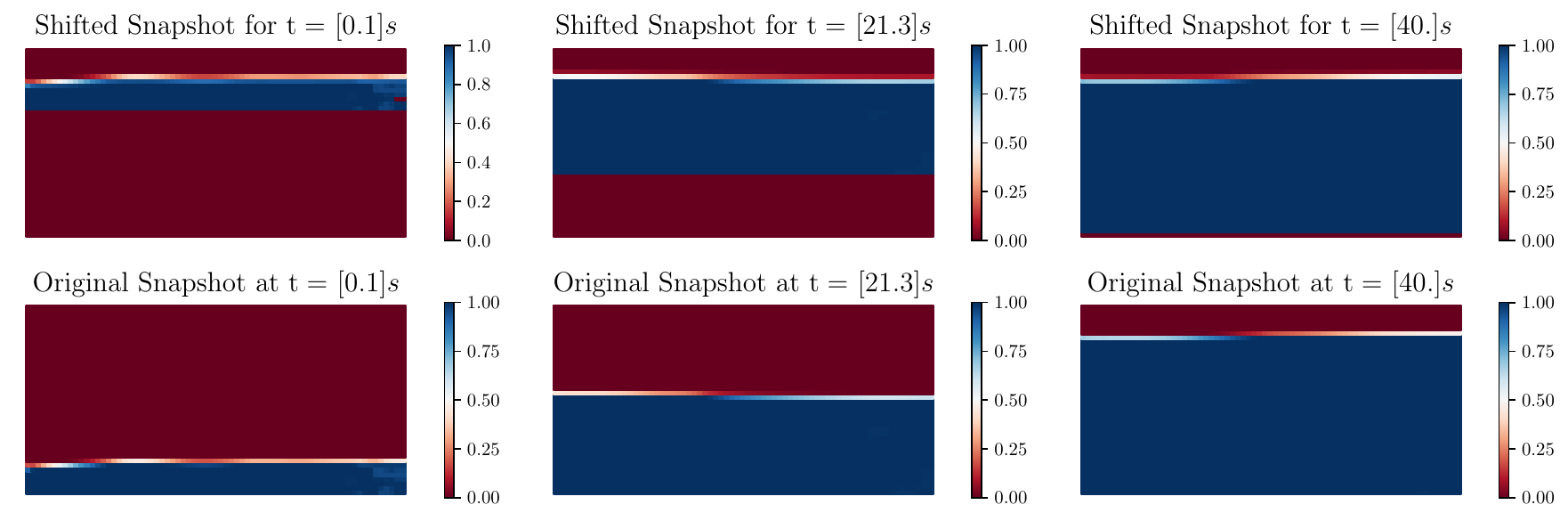}
    \centering
        \caption{Shifted snapshots (top) of two-phase flow, after neural network shift augmented manifold transformations and the respective original snapshot (bottom) at a time, t $ = \{0.05, 21.3, 40\}$ from the training set. } 
    \label{fig:two_phase_Shifted_snapshots}
\end{figure}

\begin{figure}[ht!]
    \centering
    \includegraphics[scale=0.6]{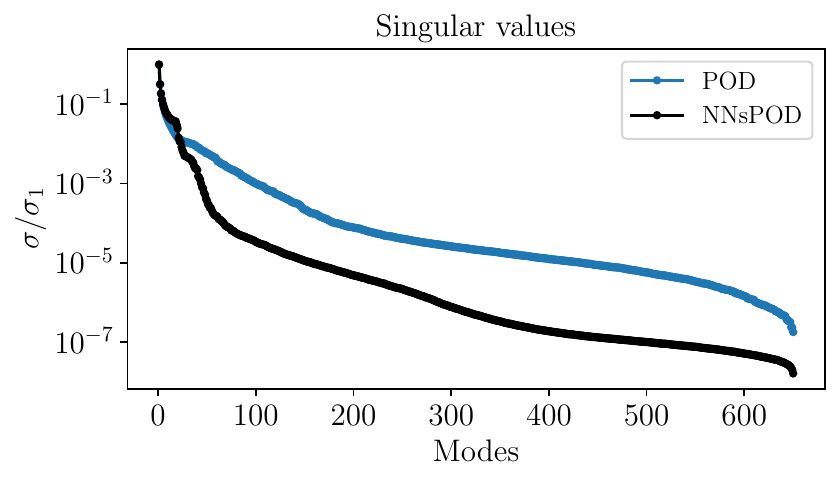}
    \caption{Comparison of singular values decay of two-phase flow FOM manifold without pre-processing, is labelled as POD and after the neural-network shifted augmented manifold transformation, is labelled as NNsPOD}
    \label{fig:two_phase_shifted_POD_singular_values}
\end{figure}

After neural-network shifted augmented manifold transformations, and performing model reduction using SVD yields singular values as shown in \autoref{fig:two_phase_shifted_POD_singular_values}. In this case, using neural network shift augmented transformations is not effective, as we don't achieve significant improvement in accelerating KnW decay. Therefore, it doesn't yield a suitable low-dimensional linear approximation subspace for efficient ROM construction. To address this, we adopt a different approach: First, we learn the shift operator and transform all the snapshots in the training set to the reference snapshot using interpolation. The customized framework for two-phase flow is presented in Algorithm \ref{alg:three}. Apart from this transformation step, the rest of Algorithm \ref{alg:one} and Algorithm \ref{alg:two} remain unchanged for the construction of efficient ROMs. \\

\RestyleAlgo{ruled}
\SetKwComment{Comment}{\(\triangleright\) }{ }
\begin{algorithm}[H]
\DontPrintSemicolon
\KwData{Inputs ($\mathbf{x}, t_i$) and corresponding snapshots $\mathbf{u}_{i} \in \textbf{X}$, where i = 0,...,$N_{\text{train}}$; reference configuration index ; InterpNet and ShiftNet architecture and respective criterion $\epsilon_{\text{interp}}, \epsilon_{\text{shift}}$ for loss functions. }
\KwResult{Transformed snapshots at reference frame and $\mathcal{T}_{\text{shift}}$ optimal shift operator. }
    
    \While{$ \mathcal{L}_{\text{InterpNet} } > \epsilon_{\text{interp}}$}{
        \textbf{InterpNet}.forward ; \\
        compute  $ \mathcal{L}_{\text{InterpNet}} $ ; \\
        \textbf{InterpNet}.backward ; \\
    }
    
    $\mathcal{T}_{\text{InterpNet}}$ = $ \textbf{InterpNet}$.forward; \Comment*[r]{reconstructs field values}
    
    \For{$\mathbf{u}(t_i)\in \mathbf{X}$, where i = 0 $\mathbf{to}$ $N_{\text{train}}$}{
        \While{$\mathcal{L}_{ShiftNet} > \epsilon_{\text{shift}}$}{
            ${\mathbf{\tilde{x}}_{i}}$ = $\mathbf{x}_{i}$ - \textbf{ShiftNet}.forward ; \Comment*[r]{x is physical space (x = x$_{\text{ref}}$)}
            $ \tilde{\mathbf{u}}_i=\mathcal{T}_{\text{interp}}({\mathbf{\tilde{x}}_{i}}); $ \Comment*[r]{reconstruct field values at each shifted space} 
            compute  $\mathcal{L}_{ShiftNet}$ ; \\
            \textbf{ShiftNet}.backward ; \\
        }
        $\mathcal{T}_{\text{shift}}$=\textbf{ShiftNet}.forward; \Comment*[r]{optimal shift operator}
        $\mathcal{L} = \mathcal{I}({\mathbf{x}}_{i},{\mathbf{u}}_{\text{ref}})$ ; \Comment*[r]{interpolator}
        $\tilde{\mathbf{u}}_{i}^{ref} = \mathcal{L}(\mathbf{x}) $; \Comment*[r]{transformed snapshots} 
    }
\caption{Customized framework for the two-phase flow case.}
\label{alg:three}
\end{algorithm}

\begin{figure}[h!]
    \includegraphics[scale=0.5]{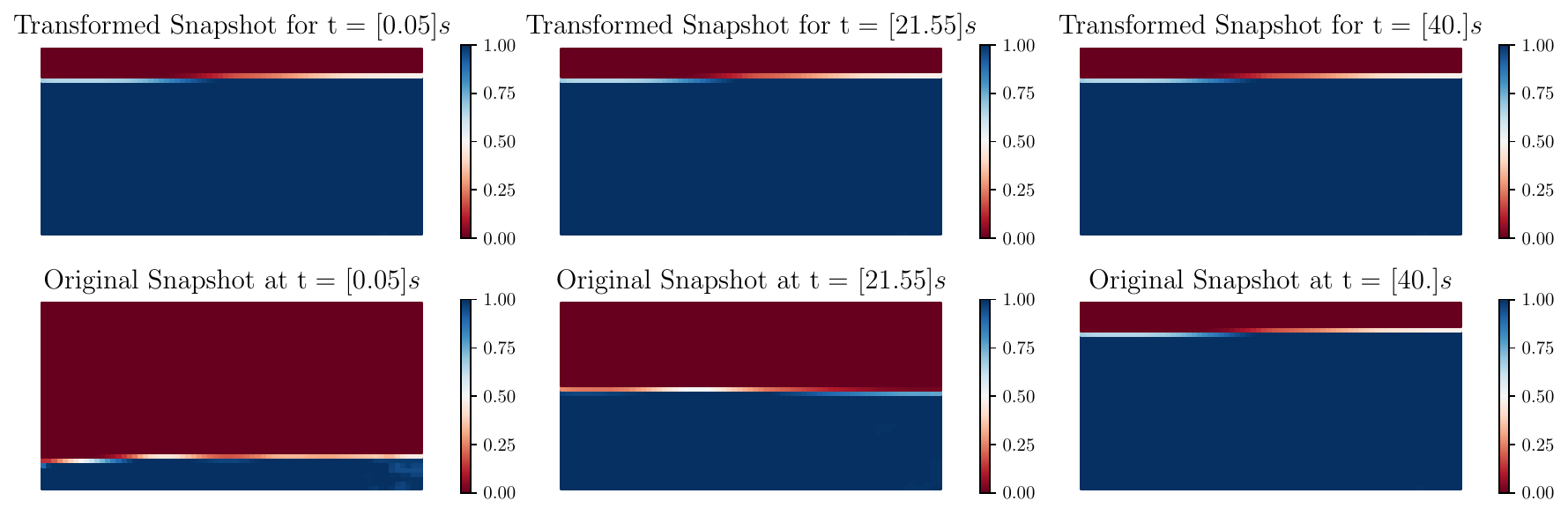}
    \centering
        \caption{Transformed snapshots (top) of two-phase flow, after transforming all the snapshots to reference snapshot using interpolation to accelerate KnW decay and the respective original snapshot (bottom) from the training set at times, t $= \{0.05, 21.55, 40\}$ are shown.} 
    \label{fig:two_phase_transformed_snapshots}
\end{figure}

\begin{figure}[h!]
    \centering
    \includegraphics[scale=0.75]{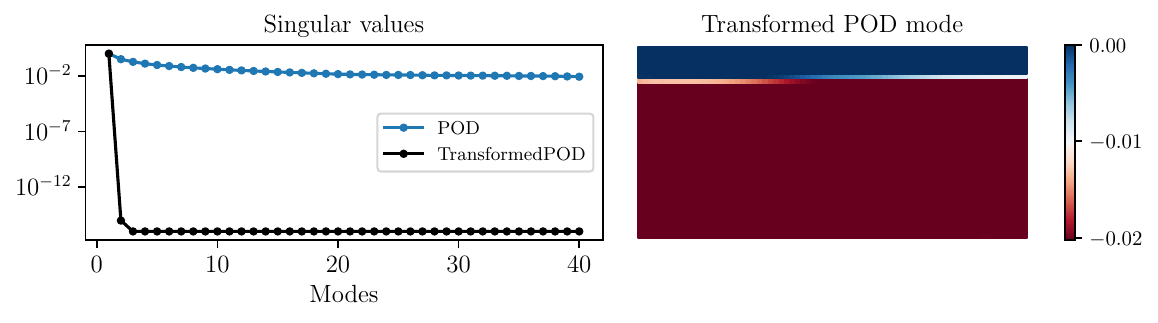}
    \caption{ \textbf{Left}: Comparison of singular values decay of two-phase flow case FOM manifold without transformations, is labelled as POD and after transforming all the snapshots to reference snapshot using interpolation, after learning shift operator, is labelled as TransformedPOD, and \textbf{Right}: First TransformedPOD mode obtained after this manifold transformation.}
    \label{fig:two_phase_sing_POD_Mode}
\end{figure}

\begin{figure}[h!]
    \centering
    \includegraphics[scale=0.5]{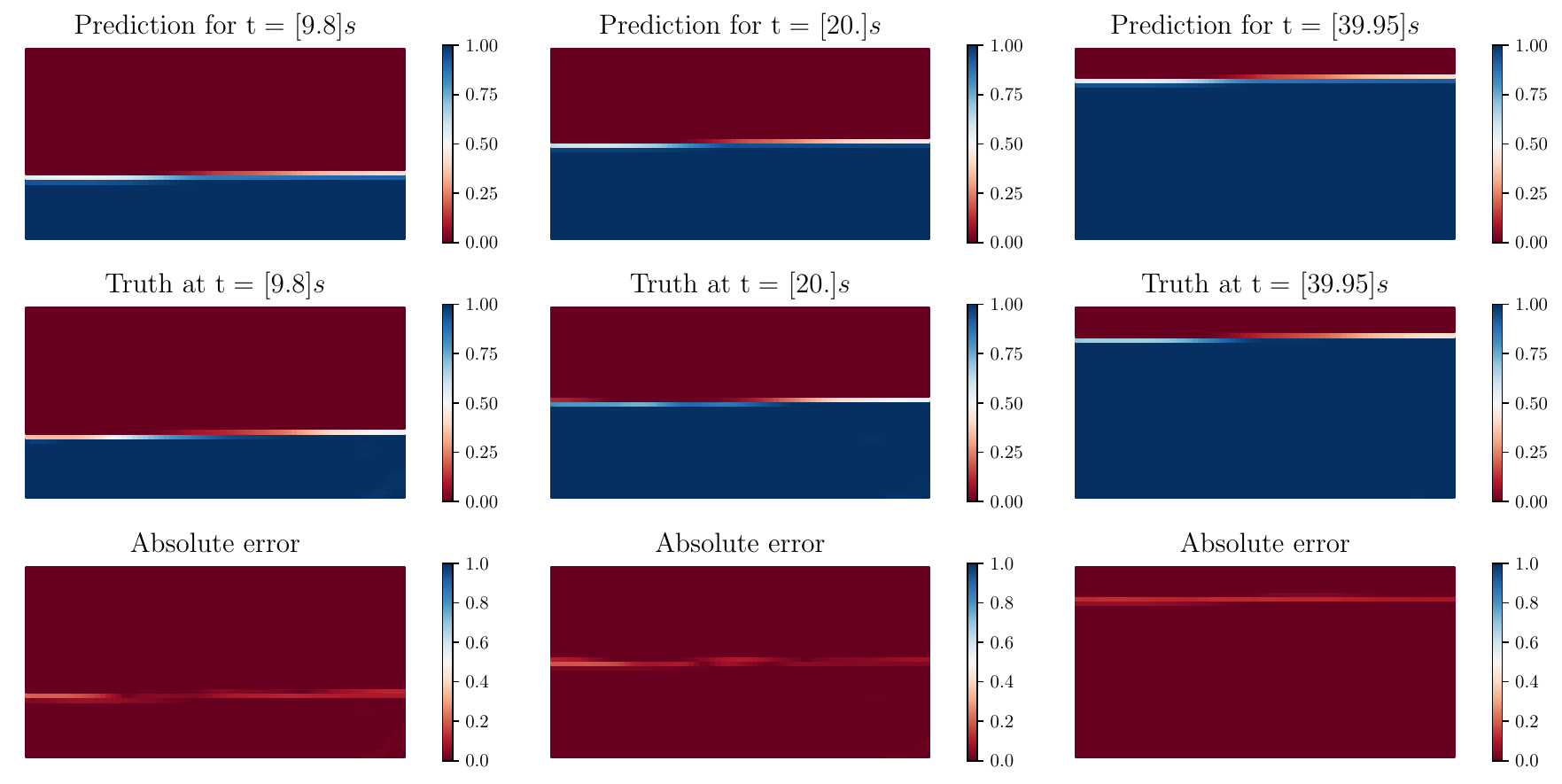}
    \caption{Predictions of snapshots for test set parameters obtained by TransformedPOD-RBF-ROM of the two-phase flow test case, showing predictions for parameters (time) t $ = \{9.8, 20, 39.95\}$s. To construct this ROM, we use the transformed linear subspace obtained by transforming all the snapshots in the training set to reference snapshot using interpolation, to accelerate KnW decay and before that we learn the shift operator. In this case, except for this transformation step in Algorithm \ref{alg:one}, the rest of Algorithm \ref{alg:one} and Algorithm \ref{alg:two} remain the same for the construction of efficient and accurate ROMs.}
    \label{fig:two_phase_prediction_test}
\end{figure}

\begin{figure}[h!]
    \centering
    \includegraphics[scale=0.53]{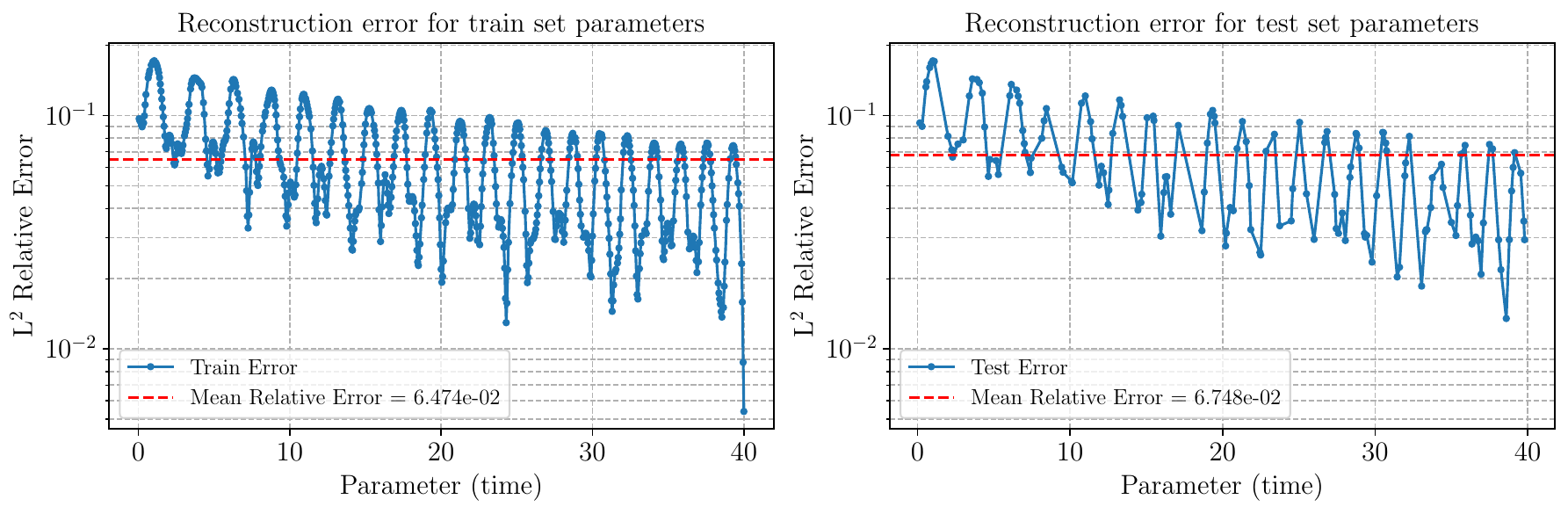}
    \caption{Relative $L^2$ reconstruction error of predictions by TransformedPOD-RBF-ROM of the two-phase flow test case, for both training set parameters (left) and test set parameters (right).} \label{fig:two_phase_prediction_error_combined}
\end{figure} 

\begin{figure}[h!]
    \centering
    \includegraphics[scale=0.5]{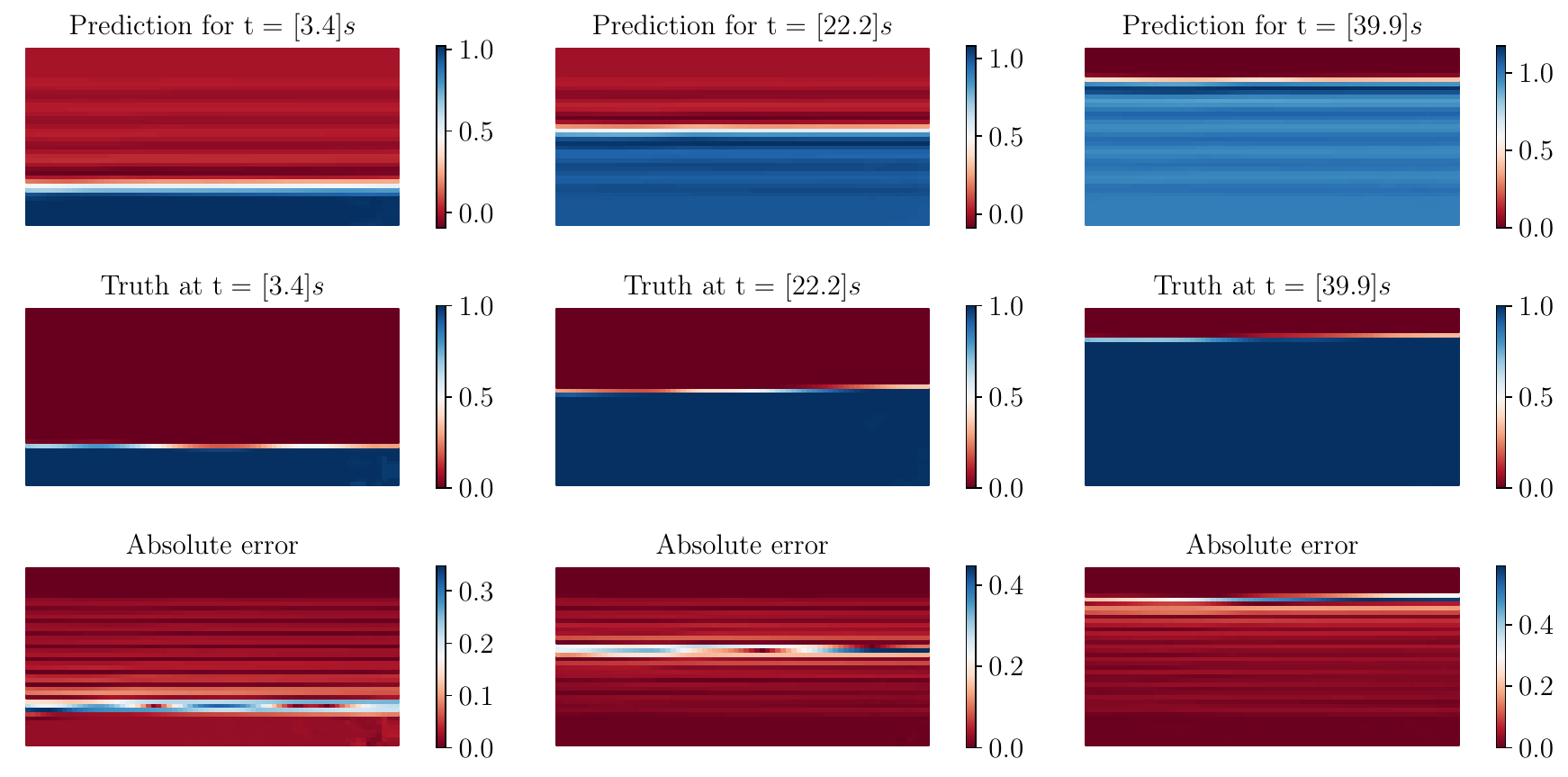}
    \caption{Predicted snapshots for test set parameters by the POD-RBF-ROM for the two-phase flow test case, considering 10 POD modes, showing predictions for parameters (time) t $ = \{3.4, 22.2, 39.9\}$s.}
    \label{fig:two_phase_prediction_test_10POD_modes}
\end{figure}

\begin{figure}[h!]
    \centering
    \includegraphics[scale=0.53]{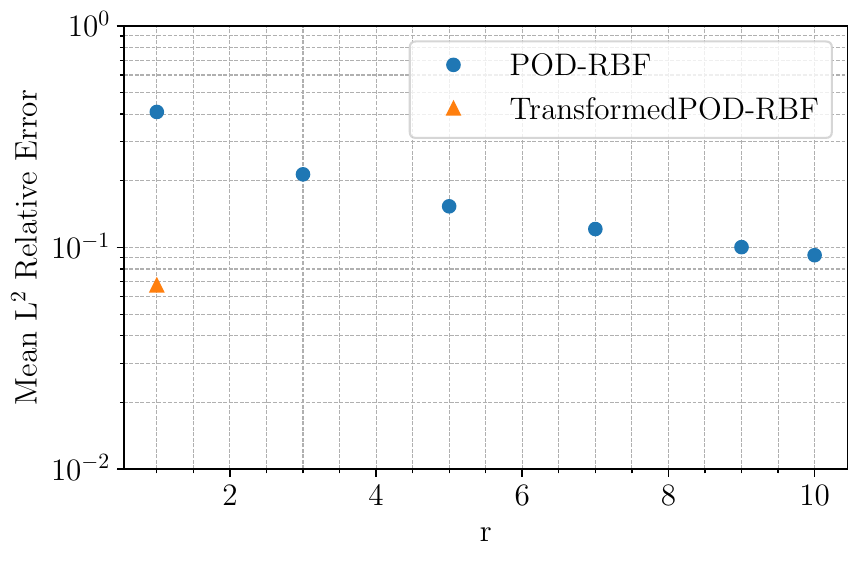}
    \caption{Mean $L^2$ relative reconstruction error of predictions by TransformedPOD-RBF-ROM of the two-phase flow test case, for both training set parameters (left) and test set parameters (right).} \label{fig:twophase_prediction_error_test_NNsPOD_POD_comparison}
\end{figure} 
\newpage
After transforming all the snapshots to reference snapshot, perform reduction using SVD, the singular values, and the first TransformedPOD mode with cumulative energy of 99.99 $\%$ is as shown in the \autoref{fig:two_phase_sing_POD_Mode}, we attain the same energy by considering 70 POD modes. For the construction of ROM, we utilize the first TransformedPOD mode, obtain the temporal coefficient, and use RBF interpolation to build, parameters to modal coefficient map. The predictions of TransformedPOD-RBF-ROM for unseen test set parameters are presented in \autoref{fig:two_phase_prediction_test}, the corresponding original snapshot (Truth), and the absolute error. It can be depicted, that the TransformedPOD-RBF-ROM predicts interface accurately and with negligible difference in the field values at the interface region. The reconstruction error for both train set and test set parameters is as shown in the \autoref{fig:two_phase_prediction_error_combined}. The error varies in the same order of magnitude for training and test set parameters. As the interface evolves over the time interval, the interface region wobbles from the left wall to the right wall of the domain (vice-versa), contributing to the observed error trend. \\

\color{black}{In \autoref{fig:two_phase_prediction_test_10POD_modes}, we illustrate the predictions for the test set parameters obtained from the POD-RBF-ROM, constructed using the first 10 POD modes. These results are qualitatively suboptimal. As discussed earlier, our objective is to develop efficient and accurate ROMs using a low-dimensional linear approximation subspace. To achieve this, we employ the first TransformedPOD mode, enabling the construction of an efficient and accurate TransformedPOD-RBF-ROM. In \autoref{fig:twophase_prediction_error_test_NNsPOD_POD_comparison}, we compare the performance of POD-RBF-ROM and TransformedPOD-RBF-ROM within the low-dimensional linear approximation subspace. While the quantitative results of the POD-RBF-ROM with 10 POD modes are comparable to those of the TransformedPOD-RBF-ROM, the latter demonstrates superior qualitative performance, highlighting its advantage over standard POD-RBF-ROMs. We highlight that the application of this method is limited to a set of problems characterized by coherent, advection-dominated features. It's not suitable for problems with significant diffusion, interface breaking, or chaotic interface dynamics. \\

We demonstrated the proposed methodology through a series of test cases with increasing complexity. In the first two cases, although the advection velocities were known, this knowledge was deliberately not utilized to perform the manifold transformation. In the two-phase flow case, where the advection velocity governing interface evolution is unknown, the proposed methodology becomes indispensable. As discussed earlier, this approach is particularly valuable for problems where determining advection velocities a priori is either impractical or impossible. It enables the derivation of a low-dimensional transformed linear approximation subspace of the transformed manifolds, efficiently approximating the FOM solutions and facilitating the construction of accurate and efficient ROMs.}
\color{black}

%% file: sections/conclusion.tex
\section{Conclusion} \label{sec:Conclusion}
The main focus of this work is the development of a complete NNsPOD-ROM algorithm for advection-dominated problems, comprising of both offline-online stages. To extend the applicability of the automatic shift detection procedure of the NNsPOD, which introduces neural-network shift augmented manifold transformations, this approach eliminates the need for comprehensive knowledge of the underlying equations and advection velocities of the given problem. It aims to derive a better linear approximation subspace of the transformed solution manifold, enabling the construction of efficient ROMs for advection-dominated problems on the transformed linear subspace. The key contribution of this work is the construction of ROM on the transformed linear subspace yielded by NNsPOD, the employment of automatic shift detection in the online stage, and we propose the complete NNsPOD-ROM algorithm for the construction of efficient ROMs. Furthermore, we tested our proposed algorithm on advection-dominated problems of increasing complexity: 1D transport equation, higher order benchmark test case - 2D isentropic convective vortex and 2D two-phase flow test cases. For 1D transport equation and 2D isentropic convective vortex test cases, the NNsPOD-ROM algorithm results in efficient and accurate ROMs. For the 2D two-phase flow test case, except for the minor adjustment in the transformation step that was required due to the specific configuration of the test case, the rest of the algorithm remains the same resulting in efficient ROM, predicts accurately the interface location, with negligible small differences in field values at the interface region. In the two-phase numerical experiment, the inlet and outlet boundaries are removed because the learned shift operator maps discrete solution points to a reference frame via spatial translation. This mapping is not a diffeomorphism and therefore does not preserve the geometry of complex domains. As a result, geometric inconsistencies may arise if such features are retained. Addressing this limitation by developing a diffeomorphic mapping, as proposed in \cite{Zhao2025}, is an important direction for future work. \\ 

Based on the numerical experiments we performed, the proposed approach yields efficient and accurate ROMs for advection-dominated flow problems. It demonstrates excellent generalization, with the training and test set errors exhibiting similar distributions across all cases. While our current experiments focus on unidirectional, the method remains effective even when the solution exhibits shifts in multiple directions. The framework is general and can be extended to handle multidirectional solutions (for example, it is demonstrated in \cite{Krah2025}). The flexibility of the approach allows for different neural network architectures to be employed depending on the specific characteristics and complexity of the problem. This flexibility could be leveraged to address such multidirectional shift scenarios, which we plan to explore in future work. Also, extending the proposed method to parametric-time dependent reduced order models is conceptually straightforward. For example, one could parameterize the inlet velocity in the two-phase flow setting. While this is not explored in the current work, we consider it a promising direction for future work. Nevertheless, the application of this approach is limited to specific advection-dominated problems, where the advected features/structure shape should be the same, and this approach is not suitable for advection-diffusion problems. The application of this approach can be extended to developing ROM for particle tracking, for which different neural-network architectures like long short-term memory (LSTM) or transformers (well suited for time series data) can be employed to learn the optimal shift operator, we postpone this for our future work. 

%% file: sections/ackno.tex
\section*{\large CRediT authorship contribution statement}
\textbf{Harshith Gowrachari}: Writing - original draft, Conceptualization, Data curation, Formal Analysis, Visualization, Methodology, Software. \textbf{Nicola Demo}: Writing – review $\&$ editing, Methodology, Software. \textbf{Giovanni Stabile}: Writing – review $\&$ editing, Methodology, Supervision. \textbf{Gianluigi Rozza}: Funding acquisition, Project administration, Supervision.  

\section*{\large Declaration of Generative AI and AI-assisted technologies in the writing process}
These technologies were used to improve readability and correct spelling during the preparation of this manuscript. After using them, the authors reviewed and edited the content as needed and take full responsibility for the content of the publication. 

\section*{\large Declaration of competing interest}
The authors declare that they have no known competing financial interests or personal relationships that could have appeared to influence the work reported in this paper

\section*{\large Acknowledgements}
We acknowledge the PhD grant supported by industrial partner Danieli \& C. S.p.A. and Programma Operativo Nazionale Ricerca e Innovazione 2014-2020, P.I. Gianluigi Rozza. GS acknowledges the financial support under the National Recovery and Resilience Plan (NRRP), Mission 4, Component 2, Investment 1.1, Call for tender No. 1409 published on 14.9.2022 by the Italian Ministry of University and Research (MUR), funded by the European Union – NextGenerationEU– Project Title ROMEU – CUP P2022FEZS3 - Grant Assignment Decree No. 1379 adopted on 01/09/2023 by the Italian Ministry of Ministry of University and Research (MUR) and acknowledges the financial support by the European Union (ERC, DANTE, GA-101115741). Views and opinions expressed are however those of the author(s) only and do not necessarily reflect those of the European Union or the European Research Council Executive Agency. Neither the European Union nor the granting authority can be held responsible for them.

\section*{\large Data statement}
Data will be made available on request. We highlight that the full order model solutions (data) were generated using open-source software, \text{HyPar} \cite{Hypar1} for the 2D Isentropic Convective Vortex, \text{OpenFoam} solver \text{interFoam} \cite{Deshpande2012} for the 2D Two-phase flow, geometry and physics inspired by industrial scenario, and the solutions of linear advection equation were obtained in closed form in \text{Python}. As reported earlier, the implementation of Algorithm \ref{alg:one} and Algorithm \ref{alg:two} is carried out using the \text{EZyRB} package (\href{https://github.com/mathLab/EZyRB}{{https://github.com/mathLab/EZyRB}}) \cite{demo18ezyrb}, an open-source \text{Python} library for non-intrusive data-driven model order reduction of parameterized problems, developed at SISSA mathLab and here, the deep neural networks are implemented using the \text{PyTorch} package \cite{paszke2017automatic}.

%% file: sections/appendix.tex
\appendix
\section{Scaling of data}\label{subsec:scaling}
Before starting the training of neural networks, we scale the data ranging between the interval [0,1] to improve the prediction accuracy of neural networks. We use Min-Max scaling for the scaling of training data between the interval [0, 1]. The Min-Max scaling for the given matrix \textit{U} is:    
\begin{equation}
U_{ij}^{\text{scaled}}=\frac{U_{ij}-\min\limits_{j=1, \ldots, N_s}\left(U_{ij}\right)}{\max\limits_{j=1, \ldots, N_s}\left(U_{ij}\right)-\min\limits_{j=1, \ldots, N_s}\left(U_{ij}\right)}, 
\end{equation}

\noindent where $i = 1, 2,......,N_h$, is the index of cell centroids or nodes in the discretized domain $\Omega$. The minimum and maximum values are stored and used to rescale the predicted results back to the original scale. 

\section{Neural networks' Architecture}\label{subsec:architecture}
In our numerical experiments, we utilize artificial neural networks, specifically the multi-layer perceptron (MLP), to implement both InterpNet and ShiftNet. The architecture and hyperparameters used for each test case are summarized below.

\begin{table}[h!]
\begin{center}
\begin{tabular}{c|cc} 
\toprule
Neural network parameters & InterpNet & ShiftNet \\
\midrule
Layers     &   [10, 10]    &   \textemdash        \\
Activation function                     &   SoftPlus                 &   LeakyReLU          \\
Learning rate                           &   0.03                     &   0.0023             \\
Accuracy threshold                      &   $10^{-6}$                &   $10^{-2}$          \\
Optimizer                               &   Adam                     &   Adam               \\
\bottomrule
\end{tabular}
\end{center}
\caption{1D Travelling wave test case.}
\label{table:1Dwave}
\end{table}

\begin{table}[h!]
\begin{center}
\begin{tabular}{c|cc} 
\toprule
Neural network parameters & InterpNet & ShiftNet \\
\midrule
Layers     &   [40, 40]    &   [10, 4]  \\
Activation function                     &   Sigmoid                     &   LeakyReLU          \\
Learning rate                           &   0.0023                      &   0.03             \\
Accuracy threshold                      &   $1.5 \times 10^{-6}$        &   0.0024          \\
Optimizer                               &   Adam                        &   Adam               \\
\bottomrule
\end{tabular}
\end{center}
\caption{2D Isentropic Convective Vortex test case.}
\label{table:Isentropic Convective Vortex}
\end{table}

\begin{table}[h!]
\begin{center}
\begin{tabular}{c|cc} 
\toprule
Neural network parameters & InterpNet & ShiftNet \\
\midrule
Layers     &   [30, 30]     &   [10, 4]  \\
Activation function                     &   LeakyReLU                   &   LeakyReLU          \\
Learning rate                           &   0.0023                      &   0.03             \\
Accuracy threshold                      &   $ 5 \times 10^{-7} $          &   0.0005         \\
Optimizer                               &   Adam                        &   Adam               \\
\bottomrule
\end{tabular}
\end{center}
\caption{2D Two-phase flow test case.}
\label{table:Two-phase}
\end{table}

\newpage